\documentclass{amsart}

\usepackage{array}
\usepackage[dvips]{color}
\usepackage{amsfonts}
\usepackage{amssymb, latexsym, amsmath, pb-diagram}
\usepackage{pstricks-add}
\usepackage{lamsarrow, pb-lams}
\usepackage{multicol}
\usepackage{bm}
\usepackage[mathscr]{eucal}
\usepackage{xy}
\usepackage{epic,eepic}
\xyoption{all}

\allowdisplaybreaks

\def\qb{\hfill $\Box$}

\newcommand\bigzero{\smash{\hbox{\large 0}}}

\begin{document}

\title{ Artin braid groups and spin structures}

\author{Gefei Wang}

\email{gefeiw@ms.u-tokyo.ac.jp}

\address{Graduate School of Mathematical Science, University of Tokyo, 3-8-1 Komaba, Meguro-ku, Tokyo, 153-8914, Japan}

\thanks{This project is supported by NSFC No.11871284}

\subjclass[2020]{20F36, 57K20, 14H30}

\keywords{spin structure, Dehn twist, $\mathfrak{S}_{2g+2}$ action, Artin braid group}

\begin{abstract}
We study the action of the Artin braid group $B_{2g+2}$ on the set of spin structures on a
hyperelliptic curve of genus $g$, which reduces to that of the symmetric group $S_{2g+2}$. It has
been already described in terms of the classical theory of Riemann surfaces. In this paper,
we compute the $S_{2g+2}$-orbits of the spin structures of genus $g$ and the isotropy group $\mathfrak{G}_i$
of each orbit in a purely combinatorial way.
\end{abstract}

\maketitle

\section{Introduction}
Let $\Sigma_g$ be a compact connected oriented surface of genus $g$, which admits a branched $2$-fold covering $\varphi: \Sigma_g \rightarrow S^2$ with $2g+2$ branch points.
If we regard $S^2$ as the Riemann sphere $\mathbb{C}P^1$, the map $\varphi$ defines the structure of a hyperelliptic curve on the surface $\Sigma_g$. The branched covering induces a group homomorphism of the Artin braid group $B_{2g+2}$ to the mapping class group of genus $g$.

A $2$-Spin structure $c$, or spin structure $c$ for short, on the surface $\Sigma_g$
can be regarded as a map $$c :  H_1(\Sigma_g , \mathbb{Z}/2) \rightarrow \mathbb{Z}/2$$ such that $c(x+y)=c(x)+c(y)+x \cdot y$ for any $x,y \in H_1(\Sigma_g , \mathbb{Z}/2)$, where $x\cdot y$ is the intersection number. The mapping class group of genus $g$ acts on the set of spin structures in a natural way, and the orbits of the action are classified by the
Arf invariant $\in \mathbb{Z}/2$.

We focus on the action of the Artin braid group of the set of spin structures of a hyperelliptic curve. By the classical theory of Riemann surfaces, the set of spin structures
is naturally bijective to the set of subsets of the branch points of the covering $\varphi$ whose cardinality is congruent to $g+1$ modulo $2$ ({\it cf.} e.g., \cite{M84} Proposition 3.1, p.3.95). In particular, the action of the Artin braid group is reduced to that of the symmetric group $\mathfrak{S}_{2g+2}$, and the isotropy group of each spin structure is isomorphic to that of the corresponding subset of the set of branch points. In particular, the action has $\lceil g/2\rceil+1$ orbits.
If $g=2k-1$, The set of spin structures under $\mathfrak{S}_{2g+2}$ action exists a unique fixed point, the isotropy group $\mathfrak{G}_i$ of the rest orbits can be described by the direct products of two symmetric groups $\mathfrak{S}_{2k+2i} \times \mathfrak{S}_{2k-2i}$ or a $\mathbb{Z}/2$-extension of direct products of two symmetric groups $\mathfrak{S}_{2k} \overleftrightarrow{\times} \mathfrak{S}_{2k}$. If $g=2k$, the isotropy group $\mathfrak{G}_i$ of each orbits can be described by the direct products of two symmetric groups $\mathfrak{S}_{2k+1+2i} \times \mathfrak{S}_{2k+1-2i}$ or a $\mathbb{Z}/2$-extension of direct products of two symmetric groups $\mathfrak{S}_{2k+1} \overleftrightarrow{\times} \mathfrak{S}_{2k+1}$.

In this paper, we will give an alternative proof of these results in a purely combinatorial way, which based on the fact that each spin structure can be regarded as a map $c$ on $H_1(\Sigma_g , \mathbb{Z}/2)$. We fix a symplectic basis
$\{
\alpha_1, \alpha_2, \dots, \alpha_g,
\beta_1, \beta_2, \dots, \beta_g\}$ for the $\mathbb{Z}/2$ first homology group $H_1(\Sigma_g, \mathbb{Z}/2)$, and introduce  a $2\times g$ $\mathbb{Z}/2$ matrix $$
M =
\left(
  \begin{array}{cccc}
    c(\alpha_1) & c(\alpha_2) & \cdots & c(\alpha_g) \\
    c(\beta_1) & c(\beta_2) & \cdots & c(\beta_g)
  \end{array}
\right).$$
In terms of the matrix, we introduce two kinds of normal forms $\{M_i\}_{i=1}^{\lceil g/2\rceil+1}$ and $\{\overline{M}_i\}_{i=1}^{\lceil g/2\rceil+1}$, where $M_i$ and $\overline{M}_i$ are in the same $\mathfrak{S}_{2g+2}$ orbit for each $i$.
In Theorem $4.2$ we prove that the any matrix $M$ is in the same $\mathfrak{S}_{2g+2}$ orbit as one of $M_i$, so that the number of $\mathfrak{S}_{2g+2}$ orbit is at most $\lceil g/2\rceil+1$. In Theorem $4.4$ and Theorem $4.5$ we prove that the isotropy group $\mathfrak{G}_i$ at $\overline{M}_i$ includes a subgroup
which can be described by the direct products of two symmetric groups or a $\mathbb{Z}/2$-extension of direct products of two symmetric groups. The normal forms $\{\overline{M}_i\}_{i=1}^{\lceil g/2\rceil+1}$ are given by Proposition $4.3$, which are different when $g=4k-1$, $g=4k$, $g=4k+1$ and $g=4k+2$, so that we need to compute the isotropy group $\mathfrak{G}_i$ at $\overline{M}_i$ for these $4$ cases to prove Theorem $4.4$ and Theorem $4.5$. Especially it is complicated to prove that the isotropy group $\mathfrak{G}_0$ at $\overline{M}_0$  includes a subgroup
which can be described by a $\mathbb{Z}/2$-extension of direct products of two symmetric groups. By using these two facts, we complete our alternative proof.

In section $2$, we recall some notions on spin structures. The action of the Artin braid group $B_{2g+2}$ on the set of spin structures is reduced to that of the symmectric group $\mathfrak{S}_{2g+2}$.
In section $3$,  we describe the $\mathfrak{S}_{2g+2}$ action explicitly, shows that the case $g=2k-1$ exists a unique fixed point and give its $\mathfrak{S}_{2g+2}$-orbits in Theorem $3.5$ and Theorem $3.6$.
In section $4$, we complete our proof of Theorem $3.5$ and Theorem $3.6$. We also give the isotropy group $\mathfrak{G}_i$ of each orbit.

 The author is grateful to Professor Nariya Kawazumi for his encouragement during this
research.

\section{Spin structure}

We take a symplectic basis
$\{
\alpha_1, \alpha_2, \dots, \alpha_g,
\beta_1, \beta_2, \dots, \beta_g\}$ for the $\mathbb{Z}/2$ first homology group $H_1(\Sigma_g, \mathbb{Z}/2)$.
By using the Poincar\'{e} Duality Theorem, one can define the intersection number $\gamma_1\cdot
\gamma_2$ by
$$
\gamma_1\cdot \gamma_2 := \langle u_1\cup u_2, [\Sigma_g]\rangle
$$
for $\gamma_1, \gamma_2 \in H_1(\Sigma_g, \mathbb{Z}/2)$.
Here $u_i \in H^1(\Sigma_g, \mathbb{Z}/2)$ is the Poincar\'{e} dual of $\gamma_i$. Then we have $\alpha_i\cdot\beta_j = \delta_{i,j}$ and
$\alpha_i\cdot\alpha_j = \beta_i\cdot\beta_j = 0$.

 The spin structure $c$ of $\Sigma_g$ with coefficient $\mathbb{Z} /2$
is well-known as a map $$c :  H_1(\Sigma_g , \mathbb{Z}/2) \rightarrow \mathbb{Z}/2$$ such that $c(x+y)=c(x)+c(y)+x \cdot y$, which is determined by
$$c(x)= \sum_{k=1}^{g} a_k b_k + f(x),$$ where $x=\sum_{k=1}^g a_k\alpha_k + b_k\beta_k$ and $f$ is a $\mathbb{Z}/2$ first cohomology class of $\Sigma_g$.
We denote by $\mathcal{C}$  the set of spin structures of $\Sigma_g$.

The Dehn twist along a simple closed curve acts on the homology group
by the automorphism
$$
T_\gamma: H_1(\Sigma_g, \mathbb{Z}/2) \to H_1(\Sigma_g, \mathbb{Z}/2)
$$
defined by $T_\gamma(x) = x+(x\cdot \gamma)\gamma$, where $\gamma$
is the $\mathbb{Z}/2$ homology class of the simple closed curve.

The $\mathbb{Z}/2$ Siegel modular group is a group of automorphisms of $H_1(\Sigma_g, \mathbb{Z}/2)$ preserving the $\mathbb{Z}/2$ intersection number, which is denoted by $Sp_{2g}(\mathbb{Z}/2)$, is generated by $T_\gamma$, where $\gamma$ runs all over $H_1(\Sigma_g,\mathbb{Z}/2)$.
The group $Sp_{2g}(\mathbb{Z}/2)$ acts on the set $\mathcal{C}$ in a natural way. In particular, $T_\gamma$ acts on a spin structure $c$ by
 $$c \circ T_\gamma (x):= c(T_\gamma(x))=c(x)+ (x \cdot \gamma)c(\gamma)+ (x \cdot \gamma)^2 .$$

There is an invariant of $c$:
 $$\text{Arf}(c):=\sum_{k=1}^g c(\alpha_k)c(\beta_k) \in \mathbb{Z}/2,$$
 which is called the Arf invariant.
The $Sp_{2g}(\mathbb{Z}/2)$-orbit of $c$ is completely determined by $\text{Arf}(c)$.

 \section{$\mathfrak{S}_{2g+2}$ action on spin structures}

 The action of the Artin braid group $B_{2g+2}$ on the cohomology group $H^1(\Sigma_{g},\mathbb{Z}/2)$ is reduced to that
 of the symmetric group $\mathfrak{S}_{2g+2}$. Hence it is similar for the action on the set $\mathcal{C}$.
More precisely,  we denote by $i : \mathfrak{S}_{2g+2} \hookrightarrow Sp_{2g}(\mathbb{Z}/2)$ the action of $\mathfrak{S}_{2g+2}$ on $H^1(\Sigma_{g},\mathbb{Z}/2)$
and take $2g+1$ transpositions $\{\sigma_i=(i,i+1)| 1 \leq i \leq 2g+1\}$ as generators,
 thus $\{\sigma_i=(i,i+1)| 1 \leq i \leq 2g+1\}$ generate $\mathfrak{S}_{2g+2}$.
Moreover we take a simple arc $l : [1,2g+2] \rightarrow S^2$ which maps each $j\in[1,2g+2]\cap\mathbb{Z}$ to a branched point of the covering $\varphi$,
and mutually disjoint simple closed curves $\overline{\beta}_i$ in $S^2$, $1\leq i\leq g$, such that $\overline{\beta}_i$ intersects the arc $l([1,2g+2])$
only at the point $l(2i-\frac{1}{2})$ transversely. We choose one of the connected component of $\varphi^{-1}(\overline{\beta}_i)$ and denote it by $\beta_{i}$,
while we denote by $\alpha_{i}$ the inverse image of $l([2i-1,2i])$. Then they form a symplectic basis
$\{
\alpha_1, \alpha_2, \dots, \alpha_g,
\beta_1, \beta_2, \dots, \beta_g\}$
for the group $H_1(\Sigma_{g},\mathbb{Z}/2)$.
Then we have
$$
 i(\sigma_1)= T_{\beta_1} \text{, } i(\sigma_{2g+1})=T_{\beta_g},
$$
 and
 $$i(\sigma_{2i})=T_{\alpha_i} \text{, } i(\sigma_{2j+1})=T_{\beta_j + \beta_{j+1}} \text{ for }1 \leq i \leq g \text{, } 1 \leq j \leq g-1.$$
As is known, the homomorphism $i$ is injective ({\it cf. e.g.,} \cite{A68}).
Hence we regard
$$\mathfrak{S}_{2g+2} \cong \langle T_\gamma | \gamma=\alpha_i,\beta_1,\beta_g, \beta_j+\beta_{j+1}, 1 \leq i \leq g , 1 \leq j \leq g-1\rangle,$$
 as a  subgroup of $Sp_{2g}(\mathbb{Z}/2)$. When $g \geq 3$, $\mathfrak{S}_{2g+2} \subsetneqq Sp_{2g}(\mathbb{Z}/2$),
while $\mathfrak{S}_6 = Sp_{4}(\mathbb{Z}/2)$ for $g=2$.

Recall that $c \in \mathcal{C}$ is determined by $$c(x)=\sum_{k=1}^ga_k b_k + f(x),$$ where
$x=\sum_{k=1}^g a_k\alpha_k+b_k\beta_k$ and $f \in H^1(\Sigma_g, \mathbb{Z}/2)$.

Since  $c\in\mathcal{C}$ is uniquely determined by the value $c(\alpha_i)$, $c(\beta_i)$, $1\leq i \leq g$, one can use a $2 \times g$ $\mathbb{Z}/2$ matrix to represent $c \in \mathcal{C}$.

{\bf Definition 3.1.}
{\it Let $c$ be a spin structure of $\Sigma_g$ determined by $$c(x)=\sum_{i=1}^ga_i b_i + f(x).$$  Then we difine the corresponding matrix
of $c$, denoted by $M$, by
$$
M :=
\left(
  \begin{array}{cccc}
    c(\alpha_1) & c(\alpha_2) & \cdots & c(\alpha_g) \\
    c(\beta_1) & c(\beta_2) & \cdots & c(\beta_g)
  \end{array}
\right)=\left(
  \begin{array}{cccc}
    f(\alpha_1) & f(\alpha_2) & \cdots & f(\alpha_g) \\
    f(\beta_1) & f(\beta_2) & \cdots & f(\beta_g)
  \end{array}
\right).$$
}

By using $M$, one can compute the $\mathfrak{S}_{2g+2}$ action on $c$ by the following lemma.

{\bf Lemma 3.2.}
{\it Let $c$ be a spin structure of $\Sigma_g$ determined by $$c(x)=\sum_{k=1}^ga_k b_k + f(x).$$
 Denote by $M$ the corresponding matrix of $c$, $$M =
\left(
  \begin{array}{cccc}
    c(\alpha_1) & c(\alpha_2) & \cdots & c(\alpha_g) \\
    c(\beta_1) & c(\beta_2) & \cdots & c(\beta_g)
  \end{array}
\right).$$
Then the $\mathfrak{S}_{2g+2}$ action on $c$ is given by:
\begin{align*}
\left(
\begin{array}{cccc}
 c(\alpha_1) & \cdots & \cdots \\
 c(\beta_1) & \cdots & \cdots
\end{array}
\right) \circ\sigma_{1}&=\left\{ \begin{array}{cccc}
\left(
  \begin{array}{cccc}
 c(\alpha_1) & \cdots & \cdots\\
 c(\beta_1) & \cdots & \cdots
  \end{array}
\right) &\text{ if }c(\beta_1)=1 \\
\left(
  \begin{array}{cccc}
    c(\alpha_1)+1 & \cdots & \cdots \\
    c(\beta_1) & \cdots & \cdots
  \end{array}
\right) &\text{ if }c(\beta_1)=0
\end{array}
\right.  \\
\left(
\begin{array}{cccc}
\cdots & \cdots & c(\alpha_g)  \\
 \cdots & \cdots & c(\beta_g)
\end{array}
\right) \circ\sigma_{2g+1}&=\left\{ \begin{array}{cccc}
\left(
  \begin{array}{cccc}
\cdots & \cdots & c(\alpha_g)  \\
 \cdots & \cdots & c(\beta_g)
  \end{array}
\right) &\text{ if }c(\beta_g)=1 \\
\left(
  \begin{array}{cccc}
\cdots & \cdots & c(\alpha_g)+1  \\
 \cdots & \cdots & c(\beta_g)
  \end{array}
\right) &\text{ if }c(\beta_g)=0
\end{array}
\right.   \\
\left(
  \begin{array}{cccc}
     \cdots & c(\alpha_i) & \cdots \\
     \cdots & c(\beta_i) & \cdots
  \end{array}
\right)\circ\sigma_{2i} & =
\left\{ \begin{array}{cccc}
\left(
  \begin{array}{cccc}
    \cdots & c(\alpha_i) & \cdots \\
     \cdots & c(\beta_i) & \cdots
  \end{array}
\right) &\text{ if }c(\alpha_i)=1 \\
\left(
  \begin{array}{cccc}
   \cdots & c(\alpha_i) & \cdots \\
     \cdots & c(\beta_i)+1 & \cdots
  \end{array}
\right) &\text{ if }c(\alpha_i)=0
\end{array}
 \right.
\end{align*}
and
\begin{align*}
&
\left(
  \begin{array}{cccc}
     \cdots & c(\alpha_j) & c(\alpha_{j+1}) & \cdots \\
     \cdots & c(\beta_j) & c(\beta_{j+1})& \cdots
  \end{array}
\right) \circ \sigma_{2j+1} \\
=&
 \left\{ \begin{array}{cccc}
\left(
  \begin{array}{cccc}
\cdots & c(\alpha_j) & c(\alpha_{j+1}) & \cdots \\
     \cdots & c(\beta_j) & c(\beta_{j+1})& \cdots
  \end{array}
\right) &\text{ if }c(\beta_j)+c(\beta_{j+1})=1 \\
\left(
  \begin{array}{cccc}
   \cdots & c(\alpha_j)+1 & c(\alpha_{j+1})+1 & \cdots \\
     \cdots & c(\beta_j) & c(\beta_{j+1})& \cdots
  \end{array}
\right) &\text{ if }c(\beta_j)+c(\beta_{j+1})=0
\end{array}
\right.
\end{align*}
for $1 \leq i \leq g$, $1 \leq j \leq g-1$.
}

{\bf Proof:}
 Firstly, we consider the $T_{\alpha_i}$ action on $c$.

 By using \[c\circ T_{\alpha_i}(x)=c(x)+(x \cdot \alpha_i)c(\alpha_i)+ (x\cdot \alpha_i)^2,\]
 where $x=\sum_{k=1}^g a_k \alpha_k + b_k \beta_k$
 and $c(x)= \sum_{k=1}^{g} a_k b_k + f(x),$
 we have $x \cdot \alpha_i =b_i$ and $c(\alpha_i)=f(\alpha_i)$.
 Since $a^2=a \in \mathbb{Z}/2$ for every $a \in \mathbb{Z}/2$, we have
\begin{align*}
 c\circ T_{\alpha_i}(x)=& \sum_{k=1}^{g} a_k b_k + f(x)+ b_i(f(\alpha_i)+1)\\
 =&\sum_{k=1}^{g} a_k b_k + \sum_{k=1}^{g}( a_k f(\alpha_k)+ b_k f(\beta_k)) + b_i(f(\alpha_i)+1)\\
 =& \sum_{k=1}^{g} a_k b_k + \sum_{k=1}^{g} a_k f(\alpha_k)+ \sum_{k\neq i} b_k f(\beta_k) + b_i (f(\beta_i)+ (f(\alpha_i)+1)),
 \end{align*}

 thus the $T_{\alpha_i}$ action on $c$ changes $c(\beta_i)$ into $c(\beta_i)+c(\alpha_i)+1$.

 Similarly, we consider the $T_{\beta_i}$ action on $c$.
\begin{align*}
 c\circ T_{\beta_i}(x)=& \sum_{k=1}^{g} a_k b_k + f(x)+ a_i(f(\beta_i)+1)\\
 =&\sum_{k=1}^{g} a_k b_k + \sum_{k=1}^{g}( a_k f(\alpha_k)+ b_k f(\beta_k)) + a_i(f(\beta_i)+1)\\
 =&\sum_{k=1}^{g} a_k b_k + \sum_{k\neq i} a_k f(\alpha_k)+\sum_{k=1}^g b_k f(\beta_k)+ a_i (f(\alpha_i) + (f(\beta_i)+1)),
 \end{align*}

 thus the $T_{\beta_i}$ action on $c$ changes $c(\alpha_i)$ into $c(\alpha_i)+c(\beta_i)+1$.

 Secondly, we consider the $T_{\beta_j+\beta_{j+1}}$ action on $c$.
By using \[c\circ T_{\beta_j+\beta_{j+1}}(x)=c(x)+(x \cdot (\beta_j+\beta_{j+1}))c(\beta_j+\beta_{j+1})+ (x\cdot (\beta_j+\beta_{j+1}))^2,\]
 we have $x \cdot (\beta_j+\beta_{j+1})=a_j + a_{j+1}$, $c(\beta_j+\beta_{j+1})=f(\beta_j)+f(\beta_{j+1})$. Then we have
 \begin{align*}
 c\circ T_{\beta_j+\beta_{j+1}}(x)=& \sum_{k=1}^{g} a_k b_k + f(x)+ (a_j + a_{j+1})(f(\beta_j)+f(\beta_{j+1})+1)\\
 =&\sum_{k=1}^{g} a_k b_k + \sum_{k=1}^{g}( a_k f(\alpha_k)+ b_k f(\beta_k)) \\
 &+ a_j  (f(\beta_j)+f(\beta_{j+1})+1)\\
 & + a_{j+1}(f(\beta_j)+f(\beta_{j+1})+1)\\
 =&\sum_{k=1}^{g} a_k b_k + \sum_{k\neq j,j+1} a_k f(\alpha_k)+ \sum_{k=1}^g b_k f(\beta_k)\\
 &+a_j (f(\alpha_j)+   (f(\beta_j)+f(\beta_{j+1})+1)) \\
 &+ a_{j+1} (f(\alpha_{j+1})+ (f(\beta_j)+f(\beta_{j+1})+1)),
 \end{align*}
 thus the $T_{\beta_j+\beta_{j+1}}$ action on $c$ changes $c(\alpha_j)$ into $c(\alpha_j)+c(\beta_j)+c(\beta_{j+1})+1$ and changes $c(\alpha_{j+1})$ into
 $c(\alpha_{j+1})+c(\beta_j)+c(\beta_{j+1})+1$.

 Finally, representing $\{T_\gamma | \gamma=\alpha_i,\beta_1,\beta_g, \beta_j+\beta_{j+1}, 1 \leq i \leq g, 1 \leq j \leq g-1\}$ by transpositions
  $\{\sigma_i=(i,i+1)| 1 \leq i \leq 2g+1\}$ and $c$ by $M$, the lemma follows. \qb

{\bf Remark:} Based on this lemma, we study the $\mathfrak{S}_{2g+2}$ action on the spin structures of $\Sigma_g$ by using matrices $M$.

 From this lemma, one can get the following theorem immediately:

{ \bf Theorem 3.3.}
{\it If and only if $g \geq 3$ is odd, there exists a unique fixed point for the $\mathfrak{S}_{2g+2}$ action on the set of spin structures of $\Sigma_g$, whose corresponding matrix is
$$\left(
  \begin{array}{cccccccccc}
   1 & 1 & 1 & 1 &\cdots & 1 & 1 & 1 & 1 \\
     1 & 0 & 1 & 0 & \cdots & 0 & 1 & 0 & 1
  \end{array}
\right).$$
Here the odd columns are $\left(
                           \begin{array}{c}
                             1 \\
                             1 \\
                           \end{array}
                         \right)$
and the even columns are $\left(
                            \begin{array}{c}
                              1 \\
                              0 \\
                            \end{array}
                          \right)$.
}

 {\bf Proof:} A spin structure $c$ is a fixed point for the $\mathfrak{S}_{2g+2}$ action if and only if $c \circ \sigma_i =c$ for every transposition $\sigma_i \in \mathfrak{S}_{2g+2}$.
 Thus the corresponding matrix $$
M =
\left(
  \begin{array}{cccc}
    c(\alpha_1) & c(\alpha_2) & \cdots & c(\alpha_g) \\
    c(\beta_1) & c(\beta_2) & \cdots & c(\beta_g)
  \end{array}
\right)$$ of $c$, satisfies $M \circ \sigma_i =M$.

 Firstly, we consider $M \circ\sigma_{2i} =M$ for $1 \leq i \leq g$. From Lemma 3.2, we get $c(\alpha_i)=1$ for $1 \leq i \leq g$.

Secondly, we consider $M \circ\sigma_1=M$, thus we get $c(\beta_1)=1$.

Thirdly, we consider $M \circ\sigma_{2j+1}=M$ for $1 \leq j \leq g-1$. From Lemma 3.2, we get $g-1$ formulas $c(\beta_j)+c(\beta_{j+1})=1$ for $1 \leq j \leq g-1$.
Since $c(\beta_1)=1$, by induction, we get
$$c(\beta_j)=\left\{
\begin{array}{cccc}
1 &\text{ if j is odd}  \\
0 &\text{ if j is even}
\end{array}\right. .$$

Finally, we consider $M \circ\sigma_{2g+1}=M$, thus we get $c(\beta_g)=1$.

In conclusion, if g is odd, we completely determined the unique fixed point $c$ by $$M=
\left(
  \begin{array}{cccccccccc}
   1 & 1 & 1 & 1 &\cdots & 1 & 1 & 1 & 1 \\
     1 & 0 & 1 & 0 & \cdots & 0 & 1 & 0 & 1
  \end{array}
\right).$$
If g is even, we have $f(\beta_g)=0$ and $f(\beta_g)=1$ at the same time, thus there does not exists any fixed point. \qb

Now we will classify $\mathfrak{S}_{2g+2}$-orbits of $\mathcal{C}$ for $g \geq 3$.

{\bf Definition 3.4.}
{\it Let $c$ be a spin structure of $\Sigma_g$ and $M$ be the corresponding matrix of $c$. Then we define a part of $M$, denoted by $A_i$, $1 \leq 2i-1 \leq g$, by
 $2\times (2i-1)$ $\mathbb{Z}/2$ matrices.
$A_1:=\left(
       \begin{array}{c}
         1 \\
         1 \\
       \end{array}
     \right)
$ and for $i \geq 2$,
\[
A_i:=\left(
       \begin{array}{cccccccc}
         1 & 1 & 1 & 1 & 1 & \cdots & 1 & 1 \\
         1 & 0 & 1 & 0 & 1 & \cdots & 0 & 1 \\
       \end{array}
     \right).
\]}

{\bf Theorem 3.5.}
{\it Let $g=2k-1$, $k \geq 2$. The number of $\mathfrak{S}_{2g+2}$-orbits of spin structures of $\Sigma_g$ is $k+1$, and
the representative of each orbit $\mathcal{C}_m$, $0 \leq m \leq k$, is given by}
\begin{align*}
M_{0}=&\left(
  \begin{array}{cccc}
   0 & 0 &\cdots & 0 \\
    0 & 0 &\cdots & 0
  \end{array}
\right)=\left(
  \bigzero
\right),\\
M_{1}= & \left( \begin{array}{cccc}
   1 & 0 &\cdots & 0 \\
   1 & 0 &\cdots & 0
  \end{array}\right)=\left(
                \begin{array}{cc}
                  A_1 & \bigzero
                \end{array}
              \right)
,\\
M_{2}= & \left(\begin{array}{cccccc}
   A_2 & \bigzero
  \end{array}
\right),\\
M_{3}= & \left( \begin{array}{cccccccc}
   A_3 & \bigzero
  \end{array}
  \right),\\
  \cdots \\
M_{k}=& A_k.
\end{align*}

{\bf Theorem 3.6.}
{\it Let $g=2k$, $k \geq 2$. The number of $\mathfrak{S}_{2g+2}$-orbits of spin structures of $\Sigma_g$ is $k+1$,
and the representative of each orbit space $\mathcal{C}_m$, $0 \leq m \leq k$, is given by}
\begin{align*}
M_{0}=&\left(
  \begin{array}{cccc}
   0 & 0 &\cdots & 0 \\
    0 & 0 &\cdots & 0
  \end{array}
\right)=\left(
  \bigzero
\right),\\
M_{1}= & \left( \begin{array}{cccc}
   1 & 0 &\cdots & 0 \\
   1 & 0 &\cdots & 0
  \end{array}
\right)=\left(
                \begin{array}{cc}
                  A_1 & \bigzero
                \end{array}
              \right),\\
M_{2}= & \left(\begin{array}{cccccc}
A_2 & \bigzero
  \end{array}
\right),\\
M_{3}= & \left( \begin{array}{cccccccc}
   A_3 & \bigzero
  \end{array}
  \right),\\
  \cdots \\
M_{k}=&\left(
  \begin{array}{cccccccccc}
   A_k & \bigzero
  \end{array}
\right).
\end{align*}

\section{Proof of Theorem 3.5 and Theorem 3.6}

In this section, we will prove that the $\mathfrak{S}_{2g+2}$ action on spin structures of $\Sigma_g$ has $k+1$ orbits for $g=2k-1$ and $g=2k$, $k \geq 2$.
We will also determine the isotropy group of each orbit.

By using Lemma 3.2, one can get following Propositon:

{\bf Proposition 4.1.}
{\it
\begin{enumerate}
\item[$(a)$] If $c(\alpha_i)=0$ for some $i$, $1 \leq i \leq g$, then the transposition $\sigma_{2i}$ changes the i-th column, thus
$$\left(\begin{array}{ccc}
\cdots & 0 & \cdots \\
\cdots & 0 & \cdots
\end{array}\right)
  \text{ and }\left(\begin{array}{ccc}
\cdots & 0 & \cdots \\
\cdots & 1 & \cdots
\end{array}\right)
 $$ are in the same orbit of $\mathcal{C}$.
 \item[$(b)$] If $c(\beta_1)=0$, then the transposition $\sigma_1$ changes the first column, thus
 $$\left(\begin{array}{ccc}
0 &  \cdots \\
0 &  \cdots
\end{array}\right)
  \text{ and } \left(\begin{array}{ccc}
 1 & \cdots \\
 0 & \cdots
\end{array}\right)
 $$ are in the same orbit of $\mathcal{C}$.
 \item[$(c_1)$]If $c(\beta_j)=c(\beta_{j+1})=\delta$, then the transposition $\sigma_{2j+1}$ changes the j-th and the (j+1)-th columns, thus
\[
\left(\begin{array}{cccc}
\cdots & 0 & 1 & \cdots \\
\cdots & \delta & \delta & \cdots
\end{array}\right)
\text{ and }\left(\begin{array}{cccc}
\cdots & 1 & 0 & \cdots \\
\cdots & \delta & \delta & \cdots
\end{array}\right)
\]
are in the same orbit of $\mathcal{C}$.
If $c(\beta_j)=c(\beta_{j+1})=\cdots=c(\beta_{j'})=\delta$, then
 \[\left(\begin{array}{cccccccccc}
\cdots &c(\alpha_{j-1}) & c(\alpha_{j}) & \cdots & c(\alpha_{j'}) & c(\alpha_{j'+1}) & \cdots \\
\cdots &c(\beta_{j-1}) & \delta & \cdots & \delta & c(\beta_{j'+1}) & \cdots
\end{array}\right)\] can be changed into
\[
\left(\begin{array}{cccccccccc}
\cdots &c(\alpha_{j-1}) & 1 & \cdots & 1 & 0 &\cdots & 0 & c(\alpha_{j'+1}) &  \cdots \\
\cdots &c(\beta_{j-1}) & \delta & \cdots & \delta & \delta &\cdots & \delta & c(\beta_{j'+1}) &  \cdots
\end{array}\right)
\]
by finite times of such transpositions.
\item[$(c_2)$]If $c(\beta_j)=c(\beta_{j+1})=\delta$, then the transposition $\sigma_{2j+1}$ changes the j-th and the (j+1)-th columns, thus
\[
\left(\begin{array}{cccc}
\cdots & 1 & 1 & \cdots \\
\cdots & \delta & \delta & \cdots
\end{array}\right)
\text{ and }\left(\begin{array}{cccc}
\cdots & 0 & 0 & \cdots \\
\cdots & \delta & \delta & \cdots
\end{array}\right)
\]
are in the same orbit of $\mathcal{C}$.
\item[$(d)$] If $c(\beta_g)=0$, then the transposition $\sigma_{2g+1}$ changes the g-th column, thus
 $$\left(\begin{array}{ccc}
  \cdots & 0\\
  \cdots & 0
\end{array}\right)
  \text{ and } \left(\begin{array}{ccc}
  \cdots & 1 \\
  \cdots & 0
\end{array}\right)
 $$ are in the same orbit of $\mathcal{C}$.
\end{enumerate}}

By using Propositon 4.1, one can find out an appropriate representative of the orbit of the corresponding matrix $M$ of $c$.

{\bf Theorem 4.2.}
{\it Let $g=2k-1$ or $g=2k$, $k \geq 2$. The number of $\mathfrak{S}_{2g+2}$-orbits of spin structures of $\Sigma_g$ is at most $k+1$.}

 In order to explain an outline of the proof of Theorem 4.2, we will show the two following examples.

{\bf Example 1} Let $g=5$ and $$M=\left(
                                    \begin{array}{ccccc}
                                      1 & 1 & 1 & 1 & 1 \\
                                      1 & 0 & 1 & 1 & 1 \\
                                    \end{array}
                                  \right)
$$ be the corresponding matrix of $c_1$.

By the $\sigma_9$ action on $M$, we get $$M \circ \sigma_{9}=M'=\left(
                                    \begin{array}{ccccc}
                                      1 & 1 & 1 & 0 & 0 \\
                                      1 & 0 & 1 & 1 & 1 \\
                                    \end{array}
                                  \right),$$ and
 by $\sigma_8\sigma_{10}$, we get $$M' \circ \sigma_8\sigma_{10}= M_2=\left(
                                    \begin{array}{ccccc}
                                      1 & 1 & 1 & 0 & 0 \\
                                      1 & 0 & 1 & 0 & 0 \\
                                    \end{array}
                                  \right).
$$
Thus $M_2$ is the appropriate representative of the orbit of $c_1$. \qb

{\bf Example 2}  Let $g=6$ and $$M=\left(
                                    \begin{array}{cccccc}
                                      1 & 1 & 1 & 1 & 1 & 1 \\
                                      1 & 0 & 1 & 1 & 0 & 1 \\
                                    \end{array}
                                  \right)
$$ be the corresponding matrix of $c_2$.

By the $\sigma_7\sigma_6\sigma_{8}$ action on $M$, we get
$$M'=\left(
                                    \begin{array}{cccccc}
                                      1 & 1 & 0 & 0 & 1 & 1 \\
                                      1 & 0 & 0 & 0 & 0 & 1 \\
                                    \end{array}
                                  \right).
$$
By $\sigma_9\sigma_7\sigma_5$, we get
$$M''=\left(
                                    \begin{array}{cccccc}
                                      1 & 0 & 0 & 0 & 0 & 1 \\
                                      1 & 0 & 0 & 0 & 0 & 1 \\
                                    \end{array}
                                  \right).
$$
By $\sigma_4\sigma_6\sigma_8\sigma_{10}\sigma_{11}\sigma_{9}\sigma_7\sigma_5$, we get
$$M'''=\left(
                                    \begin{array}{cccccc}
                                      1 & 1 & 0 & 0 & 0 & 0 \\
                                      1 & 1 & 1 & 1 & 1 & 1 \\
                                    \end{array}
                                  \right).
$$
Finally, By $\sigma_3\sigma_2\sigma_4\sigma_6\sigma_8\sigma_{10}\sigma_{12}$, we get
$$M_0=\left(
                                    \begin{array}{cccccc}
                                      0 & 0 & 0 & 0 & 0 & 0 \\
                                      0 & 0 & 0 & 0 & 0 & 0 \\
                                    \end{array}
                                  \right).
$$
Thus $M_0$ is the appropriate representative of the orbit of $c_2$. \qb

Now we prove the following theorem by using Proposition 4.1.

{\bf Proof of Theorem 4.2:} We prove this theorem by finding out an appropriate representative of the orbit of $c$, whose corresponding matrix is
$$M=\left(
  \begin{array}{cccc}
    c(\alpha_1) & c(\alpha_2) & \cdots & c(\alpha_g) \\
    c(\beta_1) & c(\beta_2) & \cdots & c(\beta_g)
  \end{array}
\right)$$.

First we assume that $M$ has at least one $\left(\begin{array}{ccc}
 1  \\
 1
\end{array}\right)$ column.

Consider the smallest $1 \leq i_1 \leq g$ such that the $i_1$-th column is $\left(\begin{array}{ccc}
 1  \\
 1
\end{array}\right)$. If $i_1 >1$, then such $M$ can be changed into a matrix whose first column is  $\left(\begin{array}{ccc}
 1  \\
 1
\end{array}\right)$, denoted by $\widetilde{M_1}$,
 by the following steps:
\begin{enumerate}
 \item From Proposition 4.1 $(a)$, $(b)$, we may assume that the first column is
$\left(\begin{array}{ccc}
0  \\
0
\end{array}\right)$ and from the second to the $(i_1-1)$-th columns are
$\left(\begin{array}{ccc}
 0  \\
 0
\end{array}\right)$ or
$\left(\begin{array}{ccc}
 1  \\
 0
\end{array}\right)$, i.e.,
\[\widetilde{M_0}:=M=\left(\begin{array}{cccccc}
0 & c(\alpha_2) & \cdots & c(\alpha_{i_1-1}) & 1 & \cdots \\
0 & 0 & \cdots& 0 & 1 & \cdots
\end{array}\right). \]
\item Find out all $2 \leq p_1 < p_2 < \cdots < p_n \leq i_1-1$ such that $c(\alpha_{p_j})=1$, $1 \leq j \leq n$.
 Making the $\sigma_{2p_1+1}\sigma_{2p_1-1}\cdots\sigma_{3}\sigma_{1}$ action on $M$ first, we get
$$\widetilde{M_0^1}=\left(\begin{array}{cccccccccc}
0 & 0 & \cdots & 0 & c(\alpha_{p_2}) & \cdots & 1 & \cdots \\
0 & 0 & \cdots & 0 & 0 & \cdots & 1 & \cdots
\end{array}\right).$$
Here the columns before the $p_2$-th column $\left(\begin{array}{ccc}
 c(\alpha_{p_2})  \\
 0
\end{array}\right)$ are $\left(\begin{array}{ccc}
 0  \\
 0
\end{array}\right)$.

Making the $\sigma_{2p_2+1}\sigma_{2p_2-1}\cdots\sigma_{3}\sigma_{1}$ action on $\widetilde{M_0^1}$, we get
$$\widetilde{M_0^2}=\left(\begin{array}{cccccccccc}
0 & 0 & \cdots & 0 & c(\alpha_{p_3}) & \cdots & 1 & \cdots \\
0 & 0 & \cdots & 0 & 0 & \cdots & 1 & \cdots
\end{array}\right).$$
Here the columns before the $p_3$-th column $\left(\begin{array}{ccc}
 c(\alpha_{p_3})  \\
 0
\end{array}\right)$ are $\left(\begin{array}{ccc}
 0  \\
 0
\end{array}\right)$.

By induction, making the $\sigma_{2p_j+1}\sigma_{2p_j-1}\cdots\sigma_{3}\sigma_{1}$ action on $\widetilde{M_0^{j-1}}$ one by one, we get
$$\widetilde{M_0^n}=\left(\begin{array}{cccccccccc}
0 & 0 & \cdots & 0 & 1 & \cdots \\
0 & 0 & \cdots & 0 & 1 & \cdots
\end{array}\right).$$
Here the columns before the $i_1$-th column $\left(\begin{array}{ccc}
 1  \\
 1
\end{array}\right)$ are $\left(\begin{array}{ccc}
 0  \\
 0
\end{array}\right)$.
\item Making the $\sigma_{2}\sigma_{4}\cdots\sigma_{2i_1-2}\cdot\sigma_{2i_1-1}\sigma_{2i_1-3}\cdots\sigma_{3}\cdot\sigma_{4}\sigma_{6}\cdots\sigma_{2i_1}$ action on $\widetilde{M_0^n}$, we get
$$\widetilde{M_1}=\left(\begin{array}{cccccccccc}
1 & 0 & \cdots & 0 & 0 & \cdots \\
1 & 0 & \cdots & 0 & 0 & \cdots
\end{array}\right).$$
Here the $i_1$-th column $\left(\begin{array}{ccc}
 1  \\
 1
\end{array}\right)$ moves to the first column and other columns do not change.
\end{enumerate}

Next find out all $1 \leq i_1 < i_2 < \cdots < i_m \leq g$ such that the $i_j$-th column is
$\left(\begin{array}{ccc}
 1  \\
 1
\end{array}\right)$.
By using the steps $(1)$, $(2)$, $(3)$ stated above, we may assume that
\[\widetilde{M_1}:=M=\left(\begin{array}{cccccccccc}
1 & c(\alpha_2) & \cdots & c(\alpha_{i_2-1}) & 1 & c(\alpha_{i_2+1}) & \cdots \\
1 & 0 & \cdots & 0 & 1 & 0 & \cdots
\end{array}\right). \]

We denote by $l_1$ the number of $\left(\begin{array}{ccc}
 1  \\
 0
\end{array}\right)$ columns between the first and the $i_2$-th column. We consider the following two cases according to the parity of $l_1$:

\begin{enumerate}
\item[$(i)$] $l_1$ is even. From Proposition 4.1 $(c_1)$, we may assume that $l_1>0$, \[\widetilde{M_1}=\left(\begin{array}{cccccccccc}
1 & 1 & \cdots & 1 & 0 &\cdots & 0 & 1 & c(\alpha_{i_2+1}) & \cdots \\
1 & 0 & \cdots & 0 & 0 &\cdots & 0 & 1 & 0 & \cdots
\end{array}\right).\] Since $l_1$ is even, from Proposition 4.1 $(c_2)$, making the $\sigma_{5}\sigma_{9}\cdots\sigma_{1+2l_1}$ action on $\widetilde{M_1}$, we get \[\left(\begin{array}{cccccccccc}
1 & 0  &\cdots & 0 & 1 & c(\alpha_{i_2+1}) & \cdots \\
1 & 0  &\cdots & 0 & 1 & 0 & \cdots
\end{array}\right).\]
Thus in this case, we only need to consider the case $l_1=0$, i.e., we have
\[\widetilde{M_1}=\left(\begin{array}{cccccccccc}
1 & 0  &\cdots & 0 & 1 & c(\alpha_{i_2+1}) & \cdots \\
1 & 0  &\cdots & 0 & 1 & 0 & \cdots
\end{array}\right).\]

Making the $\sigma_{4}\sigma_6\cdots\sigma_{2i_2-2}\cdot\sigma_{2i_2-1}\sigma_{2i_2-3}\cdots\sigma_{5}\cdot\sigma_{6}\sigma_8\cdots\sigma_{2i_2}$ action on $\widetilde{M_1}$, we get
\[\widetilde{M_1^1}=\left(\begin{array}{cccccccccc}
1 & 1 & 0  &\cdots & 0 & 0 & c(\alpha_{i_2+1}) & \cdots \\
1 & 1 & 0  &\cdots & 0 & 0 & 0 & \cdots
\end{array}\right).\]
Here the $i_2$-th column $\left(\begin{array}{ccc}
 1  \\
 1
\end{array}\right)$  moves to the second column and other columns do not change.

Now by the $\sigma_3\sigma_2\sigma_4$ action on $\widetilde{M_1^1}$, we get
\[\widetilde{M_0'}=\left(\begin{array}{cccccccccc}
0 & 0 & 0  &\cdots & 0 & 0 & c(\alpha_{i_2+1}) & \cdots \\
0 & 0 & 0  &\cdots & 0 & 0 & 0 & \cdots
\end{array}\right),\]
which is similar to $\widetilde{M_0}$.
Here we changed two $\left(\begin{array}{ccc}
 1  \\
 1
\end{array}\right)$ columns into two $\left(\begin{array}{ccc}
 0  \\
 0
\end{array}\right)$ columns, and there remains $\left(\begin{array}{ccc}
 1  \\
 1
\end{array}\right)$ in the $i_j'$-th columns for $1 < i'_1=i_3 < i'_2=i_4 < \cdots < i'_{m-2}=i_m \leq g$.

By using steps $(2)$, $(3)$, we get a new matrix from $\widetilde{M_0'}$ which is similar to $\widetilde{M_1}$ where $\left(\begin{array}{ccc}
 1  \\
 1
\end{array}\right)$ appears in the $i_j'$-th columns for $1 = i'_1 < i'_2 < \cdots < i'_{m-2} \leq g$.
\item[$(ii)$] $l_1$ is odd. From Proposition 4.1 $(c_1)$, we assume that $l_1 > 1$, \[\widetilde{M_1}=\left(\begin{array}{cccccccccc}
1 & 1 & \cdots & 1 & 0 &\cdots & 0 & 1 & c(\alpha_{i_2+1}) & \cdots \\
1 & 0 & \cdots & 0 & 0 &\cdots & 0 & 1 & 0 & \cdots
\end{array}\right).\]
Since $l_1$ is odd, from Proposition 4.1 $(c_2)$, make the $\sigma_{7}\sigma_{11}\cdots\sigma_{1+2l_1}$ action on $\widetilde{M_1}$, we get \[\left(\begin{array}{cccccccccc}
1 & 1 & 0 &\cdots & 0 & 1 & c(\alpha_{i_2+1}) & \cdots \\
1 & 0 & 0 &\cdots & 0 & 1 & 0 & \cdots
\end{array}\right).\]
Thus in this case, we only need to consider the case $l_1=1$, i.e., we have \[\widetilde{M_1}=\left(\begin{array}{cccccccccc}
1 & 1 & 0 &\cdots & 0 & 1 & c(\alpha_{i_2+1}) & \cdots \\
1 & 0 & 0 &\cdots & 0 & 1 & 0 & \cdots
\end{array}\right).\]

Now make the $\sigma_6\sigma_8\cdots\sigma_{2i_2-2}\cdot\sigma_{2i_2-1}\sigma_{2i_2-3}\cdots\sigma_7\cdot\sigma_8\sigma_{10}\cdots\sigma_{2i_2}$ action on $\widetilde{M_1}$, we get
\[\widetilde{M_2}=\left(\begin{array}{cccccccccc}
1 & 1 & 1 & 0 & \cdots & 0  & c(\alpha_{i_2+1}) & \cdots \\
1 & 0 & 1 & 0 & \cdots & 0  & 0 & \cdots
\end{array}\right).\]
Here the $i_2$-th column $\left(\begin{array}{ccc}
 1  \\
 1
\end{array}\right)$ moves to the third column and other columns do not change.
\end{enumerate}

  If $l_1$ is even, we can change $\widetilde{M_1}$ into a new matrix which is similar to $\widetilde{M_1}$ by the argument $(i)$. Then there remain $(m-2)$  $\left(\begin{array}{ccc}
 1  \\
 1
\end{array}\right)$ columns. Denote the new matrix by $\widetilde{M_1'}$. If the new $l'_1$ of $\widetilde{M_1'}$ is even, we iterate a similar process to $(i)$. If $l_1$ is odd, from $(ii)$, Proposition 4.1 $(c_1)$, we have \[\widetilde{M_2}=\left(\begin{array}{ccccccccccccc}
1 & 1 & 1 & 1 & \cdots & 1  & 0 & \cdots & 0 & 1 & c(\alpha_{i_3}+1) & \cdots \\
1 & 0 & 1 & 0 & \cdots & 0  & 0 & \cdots & 0 & 1 & 0 & \cdots
\end{array}\right).\]

We denote by $l_2$ the number of $\left(\begin{array}{ccc}
 1  \\
 0
\end{array}\right)$ columns between the third and the $i_3$-th column. We consider the following two cases according to the parity of $l_2$:
\begin{enumerate}
\item[$(iii)$] $l_2$ is even. From Proposition 4.1 $(c_1)$, $(c_2)$, we may assume that $l_2=0$ and
\[
\widetilde{M_2}=\left(\begin{array}{ccccccccccccc}
1 & 1 & 1 &    0 & \cdots & 0 & 1 & c(\alpha_{i_3}+1) & \cdots \\
1 & 0 & 1 &    0 & \cdots & 0 & 1 & 0 & \cdots
\end{array}\right).
\]

Make the $\sigma_8\sigma_{10}\cdots\sigma_{2i_3-2} \cdot \sigma_{2i_3-1}\sigma_{2i_3-3}\cdots \sigma_{9}\cdot \sigma_{10}\sigma_{12}\cdots\sigma_{2i_3}$ action on $\widetilde{M_2}$, we get
\[
\widetilde{M_2^1}=\left(\begin{array}{ccccccccccccc}
1 & 1 & 1 & 1 & 0 & \cdots & 0 & c(\alpha_{i_3}+1) & \cdots \\
1 & 0 & 1 & 1 & 0 & \cdots & 0 & 0 & \cdots
\end{array}\right).
\]
Here the $i_3$-th column $\left(\begin{array}{ccc}
 1  \\
 1
\end{array}\right)$  moves to the $4$-th column and other columns do not change.

Now make the $\sigma_7\sigma_6\sigma_8$ action on $\widetilde{M_2^1}$, we get
\[
\widetilde{M_1'}=\left(\begin{array}{ccccccccccccc}
1 & 1 & 0 & 0 & 0 & \cdots & 0 & c(\alpha_{i_3}+1) & \cdots \\
1 & 0 & 0 & 0 & 0 & \cdots & 0 & 0 & \cdots
\end{array}\right),
\]
which is similar to $\widetilde{M_1}$.
Here we change two $\left(\begin{array}{ccc}
 1  \\
 1
\end{array}\right)$ columns into two $\left(\begin{array}{ccc}
 0  \\
 0
\end{array}\right)$ columns, and there remains $\left(\begin{array}{ccc}
 1  \\
 1
\end{array}\right)$ in the $i_j'$-th columns for $1 = i'_1 < i'_2=i_4 < \cdots < i'_{m-2}=i_m \leq g$ now.
\item[$(iv)$] $l_2$ is odd. From Propositon 4.1 $(c_1)$, $(c_2)$, we may assume that $l_2=1$ and
\[
\widetilde{M_2}=\left(\begin{array}{ccccccccccccc}
1 & 1 & 1 & 1 & 0 & \cdots & 0 & 1 & c(\alpha_{i_3}+1) & \cdots \\
1 & 0 & 1 & 0 & 0 & \cdots & 0 & 1 & 0 & \cdots
\end{array}\right).
\]

Make the $\sigma_{10}\sigma_{12}\cdots \sigma_{2i_3-2}\cdot \sigma_{2i_3-1}\sigma_{2i_3-3}\cdots\sigma_{11}\cdot \sigma_{12}\sigma_{14}\cdots\sigma_{2i_3}$ action on $\widetilde{M_2}$,
we get
\[
\widetilde{M_3}=\left(\begin{array}{ccccccccccccc}
1 & 1 & 1 & 1 & 1 & 0 & \cdots & 0  & c(\alpha_{i_3}+1) & \cdots \\
1 & 0 & 1 & 0 & 1 & 0 & \cdots & 0  & 0 & \cdots
\end{array}\right).
\]
Here the $i_3$-th column $\left(\begin{array}{ccc}
 1  \\
 1
\end{array}\right)$ moves to the $5$-th column and other columns do not change.
\end{enumerate}

 If $l_2$ is even, we can change $\widetilde{M_2}$ into a new matrix which is similar to $\widetilde{M_1}$ from $(iii)$. Then there remains $(m-2)$  $\left(\begin{array}{ccc}
 1  \\
 1
\end{array}\right)$ columns. Denote by $\widetilde{M_1'}$ the new matrix. Consider the parity of the new $l'_1$ of $\widetilde{M_1'}$.  If $l'_1$ is even, then we iterate a similar process to $(i)$. If $l'_1$ is odd, then from $(ii)$, we have a new matrix which is similar to $\widetilde{M_2}$. Denote by $\widetilde{M_2'}$ the new matrix.
 Consider the parity of the new $l'_2$ of $\widetilde{M_2'}$. If $l'_2$ is even, iterate a similar process to the above. If $l_2$ is odd, from $(iv)$ and Proposition 4.1 $(c_1)$, we have
\[
\widetilde{M_3}=\left(\begin{array}{ccccccccccccccccc}
1 & 1 & 1 & 1 & 1 & 1 & \cdots & 1  & 0 & \cdots & 0 & 1 & c(\alpha_{i_4+1}) & \cdots \\
1 & 0 & 1 & 0 & 1 & 0 & \cdots & 0  & 0 & \cdots & 0 & 1 & 0 & \cdots
\end{array}\right).
\]

We proceed our argument by induction.
Denote by $l_j$ the number of $\left(\begin{array}{ccc}
 1  \\
 0
\end{array}\right)$ columns between the $i_{j}$-th and the $i_{j+1}$-th column, $1\leq j \leq m$, consider the parity of $l_j$ one by one.
Then one can
change the corresponding matrix $M$ into
\begin{equation*}
\widehat{M_{m}}:=\left(\begin{array}{ccccccccccccc}
1 & 1 & 1 & 1 & 1 & 1 & \cdots & 1 & 1  & c(\alpha_{2m}) & \cdots & c(\alpha_{g}) \\
1 & 0 & 1 & 0 & 1 & 0 & \cdots & 0 & 1  & 0                 & \cdots & 0
\end{array}\right),
\end{equation*}
where $2m-1 \leq g$, $m$ is the number of $\left(\begin{array}{ccc}
 1  \\
 1
\end{array}\right)$ columns in $\widehat{M_{m}}$.
From Proposition 4.1 $(d)$, $(2)$ and the symmetry of $\mathfrak{S}_{2g+2}$, $\widehat{M_{m}}$ can be changed into
\begin{equation*}
M_{m}=\left(\begin{array}{ccccccccccccc}
A_m & \bigzero
\end{array}\right).
\end{equation*}

Finally we consider the case that $M$ has no $\left(\begin{array}{ccc}
 1  \\
 1
\end{array}\right)$ columns. Then from $(1)$, $(2)$, $M$ can be changed into
\begin{equation*}
M_0=\left(\begin{array}{cccc}
 0 & \cdots & 0 \\
 0 & \cdots & 0
\end{array}\right)=\left(
 \bigzero
\right).
\end{equation*}

In conclusion, for every corresponding matrix $M$, one can change $M$ into one of
\begin{align*}
M_{0}=&\left(
  \bigzero
\right),\\
M_{1}= & \left( \begin{array}{cccc}
   A_1 & \bigzero
  \end{array}
\right),\\
M_{2}= & \left(\begin{array}{cccccc}
   A_2 & \bigzero
  \end{array}
\right),\\
  \cdots \\
M_{m}=&\left(
  \begin{array}{cccccccccc}
A_m & \bigzero
  \end{array}
\right),\\
\cdots
\end{align*}
We take $M_m$ as the representative of the orbit of $M$. Since $0 \leq m \leq k$ for $g=2k-1$ or $g=2k$, the theorem follows. \qb

Denote by $\mathcal{C}_m$ the $\mathfrak{S}_{2g+2}$-orbit of $M_m$ in $\mathcal{C}$, $0 \leq m \leq k$.

 We construct a subgroup of the isotropy group at $M_m$ by using the following Proposition.

{\bf Proposition 4.3.}
{\it
\begin{enumerate}
\item[$(a)$] Let $g=4k-1$. We define matrices $\overline{M_m}$, $0 \leq m \leq 2k$, by
\begin{align*}
\overline{M_0}=&
\left(\begin{array}{ccccccccccccccccc}
    A_k &
            \begin{array}{c}
              0 \\
              0 \\
            \end{array}
          & A_k
  \end{array}\right), \\
  \overline{M_1}=&
\left(\begin{array}{ccc}
A_k &
\begin{array}{ccc}
      1 & 0 & 1  \\
     0 & 1 & 0
  \end{array}
  & A_{k-1}
  \end{array}\right), \\
  \cdots \\
  \overline{M_{2i}}=&
\left(\begin{array}{ccc}
A_{k+i} &
\begin{array}{ccc}
  0  \\
  0
  \end{array} & A_{k-i}
  \end{array}\right), \\
  \overline{M_{2i+1}}=&
\left(\begin{array}{ccc}
A_{k+i} &
\begin{array}{ccc}
 1 & 0 & 1  \\
 0 & 1 & 0
  \end{array} & A_{k-i-1}
  \end{array}\right), \\
  \cdots \\
  \overline{M_{2k-1}}=&\left(
  \begin{array}{cc}
  A_{2k-1} &
  \begin{array}{cc}
    1 & 0 \\
    0 & 1
  \end{array}\end{array}\right),\\
  \overline{M_{2k}}=& A_{2k}
=M_{2k}.
\end{align*}
Then $\overline{M_m}$ is in the same orbit as $M_m$, $0 \leq m \leq 2k$.
\item[$(b)$] Let $g=4k+1$. We define matrices $\overline{M_m}$, $0 \leq m \leq 2k+1$, by
\begin{align*}
\overline{M_0}=&
\left(\begin{array}{ccc}
A_k &
\begin{array}{ccc}
     1 & 0 & 1 \\
     0 & 1 & 0
  \end{array}
  & A_k
  \end{array}\right), \\
  \overline{M_1}=&
\left(\begin{array}{ccc}
A_{k+1} &
\begin{array}{c}
  0 \\
  0
  \end{array}
  & A_k
  \end{array}\right), \\
  \cdots \\
  \overline{M_{2i}}=&
\left(\begin{array}{ccc}
A_{k+i} &
\begin{array}{ccc}
      1 & 0 & 1 \\
      0 & 1 & 0
  \end{array}
  & A_{k-i}\end{array}\right), \\
  \overline{M_{2i+1}}=&
\left(\begin{array}{ccc}
A_{k+i+1} &
\begin{array}{ccc}
     0 \\
     0
  \end{array}
  & A_{k-i}\end{array}\right),\\
  \cdots \\
  \overline{M_{2k}}=&\left(
  \begin{array}{cc}
  A_{2k} &
  \begin{array}{cc}
    1 & 0 \\
    0 & 1
  \end{array}\end{array}\right),\\
  \overline{M_{2k+1}}=& A_{2k+1}
=M_{2k+1}.
\end{align*}
Then $\overline{M_m}$ is in the same orbit as $M_m$, $0 \leq m \leq 2k+1$.
\item[$(c)$] Let $g=4k$, We define matrices $\overline{M_m}$, $0 \leq m \leq 2k$, by
\begin{align*}
\overline{M_0}=&
\left(\begin{array}{ccc}
A_k &
\begin{array}{cc}
      1 & 1  \\
      0 & 0
  \end{array}
  & A_k\end{array}\right), \\
  \overline{M_1}=&
\left(\begin{array}{cc}
     A_{k+1} & A_k
  \end{array}\right), \\
  \cdots \\
  \overline{M_{2i}}=&
\left(\begin{array}{ccc}
A_{k+i} &
\begin{array}{cc}
      1 & 1  \\
      0 & 0
  \end{array}
  & A_{k-i}\end{array}\right), \\
  \overline{M_{2i+1}}=&
\left(\begin{array}{cc}
     A_{k+i+1} & A_{k-i}
  \end{array}\right), \\
  \cdots \\
  \overline{M_{2k-1}}=&\left(
  \begin{array}{cc}
  A_{2k} &
  \begin{array}{c}
     1  \\
      1
  \end{array}\end{array}\right),\\
  \overline{M_{2k}}=&\left(\begin{array}{cc}
  A_{2k} &
  \begin{array}{cc}
     1  \\
     0
  \end{array}\end{array}
\right).
\end{align*}
Then $\overline{M_m}$ is in the same orbit as $M_m$, $0 \leq m \leq 2k$.
\item[$(d)$] Let $g=4k+2$. We define matrices $\overline{M_m}$, $0 \leq m \leq 2k+1$, by
\begin{align*}
\overline{M_0}=&
\left(\begin{array}{cc}
  A_{k+1} & A_{k+1}
  \end{array}\right), \\
  \overline{M_1}=&
\left(\begin{array}{ccc}
A_{k+1} &
\begin{array}{cc}
      1 & 1 \\
      0 & 0
  \end{array}&
  A_{k}
  \end{array}\right), \\
  \cdots \\
  \overline{M_{2i}}=&
\left(\begin{array}{cc}
    A_{k+i+1} & A_{k-i+1}
  \end{array}\right), \\
  \overline{M_{2i+1}}=&
\left(\begin{array}{ccc}
A_{k+i+1} &
\begin{array}{cc}
     1 & 1  \\
     0 & 0
  \end{array} &
  A_{k-i}\end{array}\right), \\
  \cdots \\
  \overline{M_{2k}}=&\left(
  \begin{array}{cc}
  A_{2k+1} &
  \begin{array}{c}
     1  \\
     1
  \end{array}\end{array}\right),\\
  \overline{M_{2k+1}}=&\left(
  \begin{array}{cc}
  A_{2k+1} &
  \begin{array}{c}
     1  \\
     0
  \end{array}\end{array}
\right).
\end{align*}
 Then $\overline{M_m}$ is in the same orbit as $M_m$, $0 \leq m \leq 2k+1$.
\end{enumerate}
}

{\bf Proof:} We only prove the case $(a)$, the proofs for the other cases are similar.

Consider $\overline{M_{2i}}$ first. There exists $2k$ $\left(
                                                                        \begin{array}{c}
                                                                          1 \\
                                                                          1 \\
                                                                        \end{array}
                                                                      \right)
$ columns in $\overline{M_{2i}}$.

From Proposition 4.1 $(a)$, $(c_1)$, $(c_2)$,
$\overline{M_{2i}}$ can be changed into
\[\left(\begin{array}{ccc}
A_{k+i-1}&
\begin{array}{cccccccccccccccccccc}
      1 & 0 & 0 & 0 & 1   \\
      0 & 0 & 0 & 0 & 0
  \end{array}& A_{k-i-1}\end{array}\right).
\]Here we change two $\left(\begin{array}{ccc}
 1  \\
 1
\end{array}\right)$ columns into two $\left(\begin{array}{ccc}
 0  \\
 0
\end{array}\right)$ columns, and there remain $(2k-2)$ $\left(
                                                                        \begin{array}{c}
                                                                          1 \\
                                                                          1 \\
                                                                        \end{array}
                                                                      \right)
$ columns.

From Proposition 4.1 $(a)$ $(c_1)$, $(c_2)$, it can be changed into \[\left(\begin{array}{ccc}
A_{k+i-2}&
\begin{array}{cccccccccccccccccccc}
      1 & 0 & 0 & 0 & 0 & 0 & 0 & 0 & 1   \\
      0 & 0 & 0 & 0 & 0 & 0 & 0 & 0 & 0
  \end{array} & A_{k-i-2}\end{array}\right).
\]
Here we change two $\left(\begin{array}{ccc}
 1  \\
 0
\end{array}\right)$ columns into two $\left(\begin{array}{ccc}
 0  \\
 0
\end{array}\right)$ columns, changed two $\left(\begin{array}{ccc}
 1  \\
 1
\end{array}\right)$ columns into two $\left(\begin{array}{ccc}
 0  \\
 0
\end{array}\right)$ columns, and there remain $(2k-4)$ $\left(
                                                                        \begin{array}{c}
                                                                          1 \\
                                                                          1 \\
                                                                        \end{array}
                                                                      \right)
$ columns.

By induction, we can change $\overline{M_{2i}}$ into
$M_{2i}$ by using Proposition 4.2 $(a)$, $(c_1)$, $(c_2)$ one by one. In fact there exist $(k-i)$ $\left(
                                                                        \begin{array}{c}
                                                                          1 \\
                                                                          1 \\
                                                                        \end{array}
                                                                      \right)
$ columns from the $(2k+2i)$-th column to the $g$-th column in $\overline{M_{2i}}$, then there remain $(2k-(2k-2i))=2i$ $\left(
                                                                        \begin{array}{c}
                                                                          1 \\
                                                                          1 \\
                                                                        \end{array}
                                                                      \right)
$ columns. Thus $\overline{M_{2i}}$ is in the same orbit as $M_{2i}$, $0 \leq i \leq k-1$.

Next we consider $\overline{M_{2i+1}}$. There exist $(2k-1)$ $\left(
                                                                        \begin{array}{c}
                                                                          1 \\
                                                                          1 \\
                                                                        \end{array}
                                                                      \right)
$ columns in $\overline{M_{2i+1}}$.

From Proposition 4.1 $(a)$, $(c_1)$, $(c_2)$, $\overline{M_{2i+1}}$ can be changed into
\[\left(\begin{array}{ccc}
A_{k+i} &
\begin{array}{cccccccccccccccccccc}
      0 & 0 & 0  \\
      0 & 0 & 0
  \end{array} & A_{k-i-1}\end{array}\right).\]
  Here we change two $\left(\begin{array}{ccc}
 1  \\
 0
\end{array}\right)$ columns into two $\left(\begin{array}{ccc}
 0  \\
 0
\end{array}\right)$ columns, changed $\left(
                                        \begin{array}{c}
                                          0 \\
                                          1 \\
                                        \end{array}
                                      \right)
 $ into $\left(
           \begin{array}{c}
             0 \\
             0 \\
           \end{array}
         \right)
 $.

 From Proposition 4.1 $(a)$, $(c_1)$, $(c_2)$, it can be changed into \[\left(\begin{array}{ccc}
 A_{k+i-1} &
 \begin{array}{cccccccccccccccccccc}
      1 & 0 & 0 & 0 & 0 & 0 & 1   \\
      0 & 0 & 0 & 0 & 0 & 0 & 0
  \end{array} & A_{k-i-2}\end{array}\right).\]
  Here we change two $\left(\begin{array}{ccc}
 1  \\
 1
\end{array}\right)$ columns into two $\left(\begin{array}{ccc}
 0  \\
 0
\end{array}\right)$ columns, and there remain $(2k-3)$  $\left(
                                                                        \begin{array}{c}
                                                                          1 \\
                                                                          1 \\
                                                                        \end{array}
                                                                      \right)
$ columns.

From Proposition 4.1 $(a)$, $(c_1)$, $(c_2)$, it can be changed into \[\left(\begin{array}{ccc}
A_{k+i-2} &
\begin{array}{cccccccccccccccccccc}
      1 & 0 & 0 & 0 & 0 & 0 & 0 & 0 & 0 & 0 & 1   \\
      0 & 0 & 0 & 0 & 0 & 0 & 0 & 0 & 0 & 0 & 0
  \end{array} & A_{k-i-3}\end{array}\right).\]
  Here we change two $\left(\begin{array}{ccc}
 1  \\
 0
\end{array}\right)$ columns into two $\left(\begin{array}{ccc}
 0  \\
 0
\end{array}\right)$ columns, changed two $\left(\begin{array}{ccc}
 1  \\
 1
\end{array}\right)$ columns into two $\left(\begin{array}{ccc}
 0  \\
 0
\end{array}\right)$ columns, and there remain $(2k-5)$  $\left(
                                                                        \begin{array}{c}
                                                                          1 \\
                                                                          1 \\
                                                                        \end{array}
                                                                      \right)
$ columns.

By induction, we can change $\overline{M_{2i+1}}$ into
$M_{2i+1}$ by using Proposition 4.1 $(a)$, $(c_1)$, $(c_2)$ one by one. In fact there exist $(k-i-1)$ $\left(
\begin{array}{c}
1 \\
1 \\
\end{array}
\right)
$ columns from the $(2k+2i+1)$-th column to the $g$-th column in $\overline{M_{2i+1}}$, then there remain $(2i+1)$ $\left(
\begin{array}{c}
1 \\
1 \\
\end{array}
\right)
$ columns, thus $\overline{M_{2i+1}}$ is in the same orbit as $M_{2i+1}$, $0 \leq i \leq k-1$. \qb

 Based on Proposition 4.3 $(a)$, $(b)$, we have:

 {\bf Theorem 4.4.}
 {\it Let $g=4k-1$ or $g=4k+1$, then
 $$\mathfrak{G}_m \cong \mathfrak{S}_{g+1+2m} \times \mathfrak{S}_{g+1-2m}=\langle \sigma_i | 1 \leq i \leq 2g+1, i \neq g+1+2m \rangle$$
  is a subgroup of the isotropy group at $\overline{M_m}$, where $1\leq m \leq \frac{g+1}{2}$,
  $$\mathfrak{G}_0 \cong \mathfrak{S}_{g+1} \overleftrightarrow{\times} \mathfrak{S}_{g+1}:=\langle \sigma_i, \tau | 1 \leq i \leq 2g+1, i \neq g+1 \rangle$$
   is a subgroup of the isotropy group at $\overline{M_0}$, where $\tau$ is the product of transpositions $(i,2g+3-i)$, $1 \leq i \leq g+1$, i.e., \[\tau=(1,2g+2)\cdot(2,2g+1)\cdot(3,2g)\cdots(g+1,g+2).\]}

{\bf Proof:} We will prove this theorem by finding out the generators of a subgroup of the isotropy group at $\overline{M_m}$.

Let $g=4k-1$. In
 $$\overline{M_{2j}}=
\left(\begin{array}{ccc}
A_{k+j} &
\begin{array}{ccc}
  0  \\
  0
  \end{array} & A_{k-j}
  \end{array}\right),$$
  the $(2k+2j)$-th column is $\left(
                                                            \begin{array}{c}
                                                              0 \\
                                                              0 \\
                                                            \end{array}
                                                          \right)$,
the other even columns are $\left(
                              \begin{array}{c}
                                1 \\
                                0 \\
                              \end{array}
                            \right)
$ and the odd columns are $\left(
                             \begin{array}{c}
                               1 \\
                               1 \\
                             \end{array}
                           \right)
$. Hence from Lemma 3.2, the transpositions $\sigma_i$, $1 \leq i \leq 2g+1$, $i \neq 4k+4j$, do not change $\overline{M_{2j}}$.
In
$$\overline{M_{2j+1}}=
\left(\begin{array}{ccc}
A_{k+j} &
\begin{array}{ccc}
 1 & 0 & 1  \\
 0 & 1 & 0
  \end{array} & A_{k-j-1}
  \end{array}\right),$$
 the $(2k+2j+1)$-th column is $\left(
                                                          \begin{array}{c}
                                                            0 \\
                                                            1 \\
                                                          \end{array}
                                                        \right)$,
the other odd columns are $\left(
                             \begin{array}{c}
                               1 \\
                               1 \\
                             \end{array}
                           \right)$ and the even columns are$\left(
                              \begin{array}{c}
                                1 \\
                                0 \\
                              \end{array}
                            \right)$. Hence from Lemma 3.2, the transpositions $\sigma_i$, $1 \leq i \leq 2g+1$, $i \neq 4k+4j+2$, do not change $\overline{M_{2j+1}}$.

Thus the transpositions $\{\sigma_i | 1 \leq i \leq 2g+1, i \neq 4k+2m=g+1+2m\}$ generate a subgroup of the isotropy group at $\overline{M_m}$,
$0 \leq m \leq 2k$.

Now consider the $\tau$ action on $$\overline{M_0}=
\left(\begin{array}{ccccccccccccccccc}
    A_k &
            \begin{array}{c}
              0 \\
              0 \\
            \end{array}
          & A_k
  \end{array}\right).$$ Since transpositions $$(i,2g+3-i)=\sigma_{i}\sigma_{i+1}\cdots\sigma_{2g+1-i}\sigma_{2g+2-i}\sigma_{2g+1-i}\cdots\sigma_{i+1}\sigma_i,$$ $\tau$ can be represented by the product of transpositions
\begin{align*}
\tau= &\sigma_1\sigma_2\cdots\sigma_{2g+1}\cdots\sigma_2\sigma_1 \\
& \cdot\sigma_2\sigma_3\cdots\sigma_{2g}\cdots\sigma_3\sigma_2 \\
& \cdots \\
& \cdot\sigma_{g}\sigma_{g+1}\sigma_{g+2}\sigma_{g+1}\sigma_{g}\\
& \cdot\sigma_{g+1}.
 \end{align*}

First we consider the action of $\sigma_1\sigma_2\cdots\sigma_{2g+1}\cdots\sigma_2\sigma_1$ on $\overline{M_0}$. Since the transpositions $\sigma_i$, $1 \leq i \leq 2g+1$, $i \neq g+1$, do not change $\overline{M_0}$, we have $$\overline{M_0}\circ\sigma_1\sigma_2\cdots\sigma_{2g+1}\cdots\sigma_2\sigma_1=
\overline{M_0}\circ \sigma_{g+1}\sigma_{g+2}\cdots\sigma_{2g+1}\cdots\sigma_{g+2}\sigma_{g+1}\cdots\sigma_{2}\sigma_{1}.$$
By using the relation $\sigma_{i}\sigma_{i+1}\sigma_{i}=\sigma_{i+1}\sigma_{i}\sigma_{i+1}$, we have
\begin{align*}
&\overline{M_0}\circ \sigma_{g+1}\sigma_{g+2}\cdots\sigma_{2g+1}\cdots\sigma_{g+2}\sigma_{g+1}\cdots\sigma_{2}\sigma_{1} \\
=& \overline{M_0}\circ \sigma_{2g+1}\sigma_{2g}\cdots\sigma_{g+2}\sigma_{g+1}\sigma_{g+2}\cdots\sigma_{2g}\sigma_{2g+1}\sigma_{g}\cdots\sigma_2\sigma_1 \\
=& \overline{M_0}\circ \sigma_{g+1}\cdot\sigma_{g+2}\cdots\sigma_{2g}\sigma_{2g+1}\sigma_{g}\cdots\sigma_2\sigma_1\\
=& \left(\begin{array}{ccccccccccccccccc}
     1 & 1 & 1 &  \cdots & 1 & 1 & 1 & 0 & 1 & 1 & 1 & \cdots  & 1 & 1 & 1 \\
     0 & 1 & 0 &  \cdots & 0 & 1 & 0 & 1 & 0 & 1 & 0 & \cdots  & 0 & 1 & 0
  \end{array}\right)\\
=& \left(
     \begin{array}{ccccc}
         \begin{array}{c}
           1 \\
           0 \\
         \end{array}
        & A_{k-1} &
                      \begin{array}{ccc}
                        1 & 0 & 1 \\
                        0 & 1 & 0 \\
                      \end{array}
         & A_{k-1} & \begin{array}{c}
           1 \\
           0 \\
         \end{array}
     \end{array}
   \right) \\
=& \overline{M_0^1}.
\end{align*}

Next consider $\overline{M_0^1} \circ \sigma_2\sigma_3\cdots\sigma_{2g}\cdots\sigma_3\sigma_2$,
since the transpositions $\sigma_i$, $2 \leq i \leq 2g$, $i \neq g+1$, do not change $\overline{M_0^1}$,
we have
$$
\overline{M_0^1} \circ \sigma_2\sigma_3\cdots\sigma_{2g}\cdots\sigma_3\sigma_2=\overline{M_0^1} \circ \sigma_{g+1}\sigma_{g+2}\cdots\sigma_{2g}\cdots\sigma_{g+2}\sigma_{g+1}
\cdots\sigma_3\sigma_2.
$$
By using the relation $\sigma_{i}\sigma_{i+1}\sigma_{i}=\sigma_{i+1}\sigma_{i}\sigma_{i+1}$, we have
\begin{align*}
&\overline{M_0^1}\circ \sigma_{g+1}\sigma_{g+2}\cdots\sigma_{2g}\cdots\sigma_{g+2}\sigma_{g+1}\cdots\sigma_{3}\sigma_{2} \\
=& \overline{M_0^1}\circ \sigma_{2g}\sigma_{2g-1}\cdots\sigma_{g+2}\sigma_{g+1}\sigma_{g+2}\cdots\sigma_{2g-1}\sigma_{2g}\sigma_{g}\cdots\sigma_3\sigma_2 \\
=& \overline{M_0^1}\circ \sigma_{g+1}\cdot\sigma_{g+2}\cdots\sigma_{2g-1}\sigma_{2g}\sigma_{g}\cdots\sigma_3\sigma_2\\
=& \left(\begin{array}{ccccccccccccccccc}
     0 & 1 & 1 &  \cdots & 1 & 1 & 1 & 0 & 1 & 1 & 1 & \cdots  & 1 & 1 & 0 \\
     1 & 0 & 1 &  \cdots & 1 & 0 & 1 & 0 & 1 & 0 & 1 & \cdots  & 1 & 0 & 1
  \end{array}\right)\\
=&\left(
    \begin{array}{ccccc}
      \begin{array}{cc}
        0 & 1 \\
        1 & 0
      \end{array}
       & A_{k-1} & \begin{array}{c}
                     0 \\
                     0
                   \end{array}
        & A_{k-1} & \begin{array}{cc}
        0 & 1 \\
        1 & 0
      \end{array}
    \end{array}
  \right)\\
=& \overline{M_0^2}.
\end{align*}

Consider $\overline{M_0^2} \circ \sigma_3\sigma_4\cdots\sigma_{2g-1}\cdots\sigma_4\sigma_3$,
since the transpositions $\sigma_i$,  $3 \leq i \leq 2g-1$, $i \neq g+1$, do not change $\overline{M_0^2}$,
we have
$$
\overline{M_0^2} \circ \sigma_3\sigma_4\cdots\sigma_{2g-1}\cdots\sigma_4\sigma_3=\overline{M_0^2} \circ \sigma_{g+1}\sigma_{g+2}\cdots\sigma_{2g-1}\cdots\sigma_{g+2}\sigma_{g+1}
\cdots\sigma_4\sigma_3.
$$
By using the relation $\sigma_{i}\sigma_{i+1}\sigma_{i}=\sigma_{i+1}\sigma_{i}\sigma_{i+1}$, we have
\begin{align*}
&\overline{M_0^2}\circ \sigma_{g+1}\sigma_{g+2}\cdots\sigma_{2g-1}\cdots\sigma_{g+2}\sigma_{g+1}\cdots\sigma_{4}\sigma_{3} \\
=& \overline{M_0^2}\circ \sigma_{2g-1}\sigma_{2g-2}\cdots\sigma_{g+2}\sigma_{g+1}\sigma_{g+2}\cdots\sigma_{2g-2}\sigma_{2g-1}\sigma_{g}\cdots\sigma_4\sigma_3 \\
=& \overline{M_0^2}\circ \sigma_{g+1}\cdot\sigma_{g+2}\cdots\sigma_{2g-2}\sigma_{2g-1}\sigma_{g}\cdots\sigma_4\sigma_3\\
=& \left(\begin{array}{ccccccccccccccccc}
     1 & 1 & 1 &  \cdots & 1 & 1 & 1 & 0 & 1 & 1 & 1 & \cdots  & 1 & 1 & 1 \\
     1 & 1 & 0 &  \cdots & 0 & 1 & 0 & 1 & 0 & 1 & 0 & \cdots  & 0 & 1 & 1
  \end{array}\right)\\
=&\left(
    \begin{array}{ccccc}
      A_1 & A_{k-1} & \begin{array}{ccc}
                        1 & 0 & 1 \\
                        0 & 1 & 0
                      \end{array}
       & A_{k-1} & A_1 \\
    \end{array}
  \right)\\
=& \overline{M_0^3}.
\end{align*}

It is notable that the first column and the g-th column are the same as $\overline{M_0}$'s first column and $g$-th column,
and they will not change any more by the remaining actions. Thus we proceed our argument without the first and the $g$-th column now.

Consider $\overline{M_0^3} \circ \sigma_4\sigma_5\cdots\sigma_{2g-2}\cdots\sigma_5\sigma_4$,
since the transpositions $\sigma_i$, $4 \leq i \leq 2g-2$, $i \neq g+1$, do not change $\overline{M_0^3}$,
we have
$$
\overline{M_0^3} \circ \sigma_4\sigma_5\cdots\sigma_{2g-2}\cdots\sigma_5\sigma_4=\overline{M_0^3} \circ \sigma_{g+1}\sigma_{g+2}\cdots\sigma_{2g-2}\cdots\sigma_{g+2}\sigma_{g+1}
\cdots\sigma_5\sigma_4.
$$
By using the relation $\sigma_{i}\sigma_{i+1}\sigma_{i}=\sigma_{i+1}\sigma_{i}\sigma_{i+1}$, we have
\begin{align*}
&\overline{M_0^3}\circ \sigma_{g+1}\sigma_{g+2}\cdots\sigma_{2g-2}\cdots\sigma_{g+2}\sigma_{g+1}\cdots\sigma_{5}\sigma_{4} \\
=& \overline{M_0^3}\circ \sigma_{2g-2}\sigma_{2g-3}\cdots\sigma_{g+2}\sigma_{g+1}\sigma_{g+2}\cdots\sigma_{2g-3}\sigma_{2g-2}\sigma_{g}\cdots\sigma_5\sigma_4 \\
=& \overline{M_0^3}\circ \sigma_{g+1}\cdot\sigma_{g+2}\cdots\sigma_{2g-3}\sigma_{2g-2}\sigma_{g}\cdots\sigma_5\sigma_4\\
=&   \left(\begin{array}{ccccccccccccccccc}
     1 & 0 & 1 &  \cdots & 1 & 1 & 1 & 0 & 1 & 1 & 1 & \cdots  & 1 & 0 & 1 \\
     1 & 0 & 1 &  \cdots & 1 & 0 & 1 & 0 & 1 & 0 & 1 & \cdots  & 1 & 0 & 1
  \end{array}\right)\\
=&\left(
  \begin{array}{ccccccc}
    A_1 & \begin{array}{c}
            0 \\
            0
          \end{array}
     & A_{k-1} & \begin{array}{c}
            0 \\
            0
          \end{array} & A_{k-1} & \begin{array}{c}
            0 \\
            0
          \end{array} & A_1 \\
  \end{array}
\right)\\
=& \overline{M_0^4}.
\end{align*}

Consider $\overline{M_0^4} \circ \sigma_5\sigma_6\cdots\sigma_{2g-3}\cdots\sigma_6\sigma_5$.
since the transpositions $\sigma_i$, $5 \leq i \leq 2g-3$, $i \neq g+1$, do not change $\overline{M_0^4}$,
we have
$$
\overline{M_0^4} \circ \sigma_5\sigma_6\cdots\sigma_{2g-3}\cdots\sigma_6\sigma_5=\overline{M_0^4} \circ \sigma_{g+1}\sigma_{g+2}\cdots\sigma_{2g-3}\cdots\sigma_{g+2}\sigma_{g+1}
\cdots\sigma_6\sigma_5.
$$
By using the relation $\sigma_{i}\sigma_{i+1}\sigma_{i}=\sigma_{i+1}\sigma_{i}\sigma_{i+1}$, we have
\begin{align*}
&\overline{M_0^4}\circ \sigma_{g+1}\sigma_{g+2}\cdots\sigma_{2g-3}\cdots\sigma_{g+2}\sigma_{g+1}\cdots\sigma_{6}\sigma_{5} \\
=& \overline{M_0^4}\circ \sigma_{2g-3}\sigma_{2g-4}\cdots\sigma_{g+2}\sigma_{g+1}\sigma_{g+2}\cdots\sigma_{2g-4}\sigma_{2g-3}\sigma_{g}\cdots\sigma_6\sigma_5 \\
=& \overline{M_0^4}\circ \sigma_{g+1}\cdot\sigma_{g+2}\cdots\sigma_{2g-4}\sigma_{2g-3}\sigma_{g}\cdots\sigma_6\sigma_5\\
=&   \left(\begin{array}{ccccccccccccccccc}
     1 & 1 & 1 &  \cdots & 1 & 1 & 1 & 0 & 1 & 1 & 1 & \cdots  & 1 & 1 & 1 \\
     1 & 0 & 0 &  \cdots & 0 & 1 & 0 & 1 & 0 & 1 & 0 & \cdots  & 0 & 0 & 1
  \end{array}\right)\\
=&\left(
    \begin{array}{ccccccc}
      A_1 & \begin{array}{cc}
              1 & 1 \\
              0 & 0
            \end{array}
       & A_{k-2} & \begin{array}{ccc}
                      1 & 0 & 1 \\
                      0 & 1 & 0
                    \end{array}
       & A_{k-2} & \begin{array}{cc}
              1 & 1 \\
              0 & 0
            \end{array} & A_1 \\
    \end{array}
  \right)\\
=& \overline{M_0^5}\\
&\cdots
\end{align*}

It is notable that the second column and the $(g-1)$-th column are the same as $\overline{M_0}$'s second column and $(g-1)$-th column,
and they will not change any more by the remaining actions. Thus we proceed our argument without the first, the second, the $(g-1)$-th, the $g$-th column now.

By induction, the $\sigma_{i}\sigma_{i+1}\cdots\sigma_{2g+1-i}\sigma_{2g+2-i}\sigma_{2g+1-i}\cdots\sigma_{i+1}\sigma_i$ action on $\overline{M_0^{i-1}}$ has
\begin{align*}
& \overline{M_0^{i-1}}\circ \sigma_{i}\sigma_{i+1}\cdots\sigma_{2g+1-i}\sigma_{2g+2-i}\sigma_{2g+1-i}\cdots\sigma_{i+1}\sigma_i\\
=&
\overline{M_0^{i-1}} \circ \sigma_{g+1}\cdot\sigma_{g+2}\cdots\sigma_{2g+1-i}\sigma_{2g+2-i}\sigma_g\cdots\sigma_{i+1}\sigma_i \\
=& \overline{M_0^i}.
\end{align*}
For $i=2j+1$ or $i=2j+2$, $1 \leq i \leq g$, the columns except the columns from the $(j+1)$-th column to the $(g-j)$-th column are the
same as those $\overline{M_0}$'s columns.

For $i=2j+1$, the $2k$-th column is
$\left(
   \begin{array}{c}
     0 \\
     1 \\
   \end{array}
 \right)
$.
\begin{enumerate}
\item If $j$ is odd, the other even columns are
$\left(
  \begin{array}{c}
    1 \\
    1 \\
  \end{array}
\right)$, the other odd columns are
$\left(
  \begin{array}{c}
    1 \\
    0 \\
  \end{array}
\right)
$.
\item If $j$ is even, the $(j+1)$-th column and the $(g-j)$-th column are
$\left(
   \begin{array}{c}
     1 \\
     0 \\
   \end{array}
 \right)
$, the other even columns are
$\left(
  \begin{array}{c}
    1 \\
    0 \\
  \end{array}
\right)$, the other odd columns are
$\left(
  \begin{array}{c}
    1 \\
    1 \\
  \end{array}
\right)
$.
\end{enumerate}

For $i=2j+2$, the $2k$-th column is
$\left(
   \begin{array}{c}
     0 \\
     0 \\
   \end{array}
 \right)
$.
\begin{enumerate}
\item[(3)] If $j$ is odd,the $(j+1)$-th column and the $(g-j)$-th column are
$\left(
   \begin{array}{c}
     0 \\
     0 \\
   \end{array}
 \right)
$, the other even columns are
$\left(
  \begin{array}{c}
    1 \\
    0 \\
  \end{array}
\right)$, the other odd columns are
$\left(
  \begin{array}{c}
    1 \\
    1 \\
  \end{array}
\right)
$.
\item[(4)] If $j$ is even, the $(j+1)$-th column and the $(g-j)$-th column are
$\left(
   \begin{array}{c}
     0 \\
     1 \\
   \end{array}
 \right)
$, the other even columns are
$\left(
  \begin{array}{c}
    1 \\
    0 \\
  \end{array}
\right)$, the other odd columns are
$\left(
  \begin{array}{c}
    1 \\
    1 \\
  \end{array}
\right)
$.
\end{enumerate}

Since $g=4k-1$, we have
\begin{align*}
&\overline{M_0^{g-1}}\circ \sigma_{g}\sigma_{g+1}\sigma_{g+2}\sigma_{g+1}\sigma_g \\
=&\overline{M_0^{g-1}}\circ \sigma_{g+1}\cdot\sigma_{g+2}\sigma_{g} \\
=&\left(\begin{array}{ccccccccccccccccc}
    A_k &
            \begin{array}{c}
              0 \\
              1 \\
            \end{array}
          & A_k
  \end{array}\right) \\
=& \overline{M_0^g}.
\end{align*}
Then we have
\begin{align*}
 \overline{M_0}\circ \tau
=&\overline{M_0^g}\circ \sigma_{g+1}\\
=&
\left(\begin{array}{ccccccccccccccccc}
    A_k &
            \begin{array}{c}
              0 \\
              0 \\
            \end{array}
          & A_k
  \end{array}\right) \\
=&\overline{M_0}.
\end{align*}

Now let $g=4k+1$. Similar to the case of $g=4k-1$, the transpositions $\{\sigma_i | 1 \leq i \leq 2g+1, i \neq 4k+2+2m=g+1+2m\}$ generate a subgroup of the isotropy group at $\overline{M_m}$,
$0 \leq m \leq 2k+1$.

We consider the $\tau$ action on
\[\overline{M_0}=
\left(\begin{array}{ccc}
A_k &
\begin{array}{ccc}
     1 & 0 & 1 \\
     0 & 1 & 0
  \end{array}
  & A_k
  \end{array}\right).
\]

Since
\begin{align*}
\tau= &\sigma_1\sigma_2\cdots\sigma_{2g+1}\cdots\sigma_2\sigma_1 \\
& \cdot\sigma_2\sigma_3\cdots\sigma_{2g}\cdots\sigma_3\sigma_2 \\
& \cdots \\
& \cdot\sigma_{g}\sigma_{g+1}\sigma_{g+2}\sigma_{g+1}\sigma_{g}\\
& \cdot\sigma_{g+1},
\end{align*}
we have
\begin{align*}
&\overline{M_0}\circ\sigma_1\sigma_2\cdots\sigma_{2g+1}\cdots\sigma_2\sigma_1\\
=&\overline{M_0}\circ \sigma_{g+1}\sigma_{g+2}\cdots\sigma_{2g+1}\cdots\sigma_{g+2}\sigma_{g+1}\cdots\sigma_{2}\sigma_{1} \\
=& \overline{M_0}\circ \sigma_{2g+1}\sigma_{2g}\cdots\sigma_{g+2}\sigma_{g+1}\sigma_{g+2}\cdots\sigma_{2g}\sigma_{2g+1}\sigma_{g}\cdots\sigma_2\sigma_1 \\
=& \overline{M_0}\circ \sigma_{g+1}\cdot\sigma_{g+2}\cdots\sigma_{2g}\sigma_{2g+1}\sigma_{g}\cdots\sigma_2\sigma_1\\
=& \left(\begin{array}{ccccccccccccccccc}
     1 & 1 & 1 &  \cdots & 1 & 1 & 1 & 0 & 1 & 1 & 1 & \cdots  & 1 & 1 & 1 \\
     0 & 1 & 0 &  \cdots & 1 & 0 & 1 & 0 & 1 & 0 & 1 & \cdots  & 0 & 1 & 0
  \end{array}\right)\\
=& \left(\begin{array}{ccccc}\begin{array}{c}
     1  \\
     0
  \end{array} & A_{k} & \begin{array}{c}
     0  \\
     0
  \end{array} & A_{k} & \begin{array}{c}
     1  \\
     0
  \end{array}\end{array}\right)\\
=& \overline{M_0^{1}}.
\end{align*}

Consider $\overline{M_0^{1}} \circ \sigma_2\sigma_{3}\cdots\sigma_{2g}\cdots\sigma_{3}\sigma_2$, we have
\begin{align*}
&\overline{M_0^{1}}\circ\sigma_2\sigma_{3}\cdots\sigma_{2g}\cdots\sigma_{3}\sigma_2\\
=&\overline{M_0^{1}}\circ \sigma_{g+1}\sigma_{g+2}\cdots\sigma_{2g}\cdots\sigma_{g+2}\sigma_{g+1}\sigma_g\cdots\sigma_{2} \\
=& \overline{M_0^{1}}\circ \sigma_{2g}\cdots\sigma_{g+2}\sigma_{g+1}\sigma_{g+2}\cdots\sigma_{2g}\sigma_{2g+1}\sigma_{g}\cdots\sigma_2 \\
=& \overline{M_0^{1}}\circ \sigma_{g+1}\cdot\sigma_{g+2}\cdots\sigma_{2g}\sigma_{g}\cdots\sigma_2\\
=& \left(\begin{array}{ccccccccccccccccc}
     0 & 1 & 1 &  \cdots & 1 & 1 & 1 & 0 & 1 & 1 & 1 & \cdots  & 1 & 1 & 0 \\
     1 & 0 & 1 &  \cdots & 0 & 1 & 0 & 1 & 0 & 1 & 0 & \cdots  & 1 & 0 & 1
  \end{array}\right)\\
=& \left(\begin{array}{ccccc}\begin{array}{cc}
     0 & 1 \\
     1 & 0
     \end{array} & A_{k-1} & \begin{array}{ccc}
                               1 & 0 & 1 \\
                               0 & 1 & 0
                             \end{array} & A_{k-1} & \begin{array}{cc}
     1 & 0 \\
     0 & 1
     \end{array}
  \end{array}\right)\\
=& \overline{M_0^{2}}.
\end{align*}

Consider $\overline{M_0^{2}} \circ \sigma_3\sigma_4\cdots\sigma_{2g-1}\cdots\sigma_4\sigma_3$, we have
\begin{align*}
&\overline{M_0^{2}}\circ\sigma_3\sigma_4\cdots\sigma_{2g-1}\cdots\sigma_4\sigma_3\\
=& \left(\begin{array}{cccccccccccccccccccccc}
     1 & 1 & 1 &  \cdots & 1 & 1 & 1 & 0 & 1 & 1 & 1 & \cdots  & 1 & 1 & 1 \\
     1 & 1 & 0 &  \cdots & 1 & 0 & 1 & 0 & 1 & 0 & 1 & \cdots  & 0 & 1 & 1
  \end{array}\right)\\
=& \left(\begin{array}{ccccc}
A_1 & A_{k} & \begin{array}{c}
     0  \\
     0
  \end{array} & A_{k} & A_1 \end{array}\right)\\
=& \overline{M_0^{3}}.
\end{align*}

Consider $\overline{M_0^{3}} \circ \sigma_4\sigma_5\cdots\sigma_{2g-2}\cdots\sigma_5\sigma_4$, we have
\begin{align*}
&\overline{M_0^{3}} \circ \sigma_4\sigma_5\cdots\sigma_{2g-2}\cdots\sigma_5\sigma_4\\
=& \left(\begin{array}{ccccccccccccccccc}
     1 & 0 & 1 &  \cdots & 1 & 1 & 1 & 0 & 1 & 1 & 1 & \cdots  & 1 & 0 & 1 \\
     1 & 0 & 1 &  \cdots & 0 & 1 & 0 & 1 & 0 & 1 & 0 & \cdots  & 1 & 0 & 1
  \end{array}\right)\\
=& \left(\begin{array}{ccccccc}
A_1 & \begin{array}{c}
        0 \\
        0
      \end{array} & A_{k-1} & \begin{array}{ccc}
                                1 & 0 & 1 \\
                                0 & 1 & 0
                              \end{array} & A_{k-1} & \begin{array}{c}
                                                        0 \\
                                                        0
                                                      \end{array} & A_1
\end{array}\right)\\
=& \overline{M_0^{4}}\\
&\cdots
\end{align*}

By induction, we make the same changes by the actions as the case $g=4k-1$. In fact, we think of the middle columns $\left(
                                              \begin{array}{ccc}
                                                1 & 0 & 1 \\
                                                0 & 1 & 0 \\
                                              \end{array}
                                            \right)
$ in $\overline{M_0^{2i}}$ as $\left(
        \begin{array}{c}
          0 \\
          0 \\
        \end{array}
      \right)
$ and think of the middle columns $\left(
                                              \begin{array}{ccc}
                                                1 & 0 & 1 \\
                                                1 & 0 & 1 \\
                                              \end{array}
                                            \right)
$ in $\overline{M_0^{2i-1}}$ as $\left(
                                              \begin{array}{ccc}
                                                 0  \\
                                                 1  \\
                                              \end{array}
                                            \right)
$.

Thus we have
\begin{align*}
&\overline{M_0^{g-1}}\circ \sigma_{g}\sigma_{g+1}\sigma_{g+2}\sigma_{g+1}\sigma_g \\
=&\overline{M_0^{g-1}}\circ \sigma_{g+1}\cdot\sigma_{g+2}\sigma_{g} \\
=&\left(
                                   \begin{array}{ccc}
                                   A_k &\begin{array}{ccc}
                                    1 & 0 & 1 \\
                                    0 & 0 & 0 \\
                                   \end{array}
                                 & A_k \end{array}\right) \\
=& \overline{M_0^g},
\end{align*}
 then we get
\begin{align*}
 \overline{M_0}\circ \tau =& \overline{M_0^{g}}\circ \sigma_{g+1}\\
 = &\left(
                                   \begin{array}{ccc}
                                   A_k &\begin{array}{ccc}
                                    1 & 0 & 1 \\
                                    0 & 1 & 0 \\
                                   \end{array}
                                 & A_k
                                 \end{array}\right)\\
 =& \overline{M_0}.
 \end{align*}

 The theorem follows. \qb

 Now consider the cases $g=4k$ and $g=4k+2$, we have the similar results.

 {\bf Theorem 4.5.}
 {\it Let $g=4k$ or $g=4k+2$, then
 $$\mathfrak{G}_m \cong \mathfrak{S}_{g+1+2m} \times \mathfrak{S}_{g+1-2m}=\langle \sigma_i | 1 \leq i \leq 2g+1, i \neq g+1+2m \rangle$$
  is a subgroup of the isotropy group at $\overline{M_m}$, where $1\leq m \leq \frac{g}{2}$,
  $$\mathfrak{G}_0 \cong \mathfrak{S}_{g+1} \overleftrightarrow{\times} \mathfrak{S}_{g+1}:=\langle \sigma_i, \tau | 1 \leq i \leq 2g+1, i \neq g+1 \rangle$$
   is a subgroup of the isotropy group at $\overline{M_0}$, where $\tau$ is the product of transpositions $(i,2g+3-i)$, $1 \leq i \leq g+1$, i.e., \[\tau=(1,2g+2)\cdot(2,2g+1)\cdot(3,2g)\cdots(g+1,g+2).\]}

{\bf Proof:} The proof is similar to the proof of Theorem 4.4.

Let $g=4k$. It is easy to see that the transpositions $\{\sigma_i | 1 \leq i \leq 2g+1, i \neq 4k+1+2m=g+1+2m\}$ generate a subgroup of the isotropy group at $\overline{M_m}$,
$0 \leq m \leq 2k$.

We consider the $\tau$ action on $$\overline{M_0}=\left(\begin{array}{ccc}
A_k &
\begin{array}{cc}
      1 & 1  \\
      0 & 0
  \end{array}
  & A_k\end{array}\right).$$
Since we have $$(i,2g+3-i)=\sigma_{i}\sigma_{i+1}\cdots\sigma_{2g+1-i}\sigma_{2g+2-i}\sigma_{2g+1-i}\cdots\sigma_{i+1}\sigma_i,$$ $\tau$ can be represented by the product of transpositions
\begin{align*}
\tau= &\sigma_1\sigma_2\cdots\sigma_{2g+1}\cdots\sigma_2\sigma_1 \\
& \cdot\sigma_2\sigma_3\cdots\sigma_{2g}\cdots\sigma_3\sigma_2 \\
& \cdots \\
& \cdot\sigma_{g}\sigma_{g+1}\sigma_{g+2}\sigma_{g+1}\sigma_{g}\\
& \cdot\sigma_{g+1}.
 \end{align*}

 First we consider the $\sigma_1\sigma_2\cdots\sigma_{2g+1}\cdots\sigma_2\sigma_1$ action on $\overline{M_0}$,
\begin{align*}
&\overline{M_0}\circ \sigma_1\sigma_2\cdots\sigma_{2g+1}\cdots\sigma_2\sigma_1 \\
=&\overline{M_0}\circ \sigma_{g+1}\cdot\sigma_{g+2}\sigma_{g+3}\cdots\sigma_{2g+1}\sigma_{g}\sigma_{g-1}\cdots\sigma_1 \\
=&\left(\begin{array}{cccc}
\begin{array}{c}
  1 \\
  0
\end{array}
 & A_k  & A_k &
\begin{array}{c}
  1 \\
  0
\end{array}
\end{array}\right) \\
=& \overline{M_0^1}.
\end{align*}

Second we consider the $\sigma_2\sigma_3\cdots\sigma_{2g}\cdots\sigma_3\sigma_2$ action on $\overline{M_0^1}$,
\begin{align*}
&\overline{M_0^1}\circ \sigma_2\sigma_3\cdots\sigma_{2g}\cdots\sigma_3\sigma_2 \\
=&\overline{M_0^1}\circ \sigma_{g+1}\cdot\sigma_{g+2}\sigma_{g+3}\cdots\sigma_{2g}\sigma_{g}\sigma_{g-1}\cdots\sigma_2 \\
=&\left(\begin{array}{ccccc}
\begin{array}{cc}
  0 & 1 \\
  1 & 0
\end{array}
 & A_{k-2}  & \begin{array}{cc}
                1 & 1 \\
                0 & 0
              \end{array}
 & A_{k-2} &
\begin{array}{cc}
 0 & 1 \\
 1 & 0
\end{array}
\end{array}\right) \\
=& \overline{M_0^2} \\
& \cdots
\end{align*}

{ In fact, if we add  a  $\left(
                                                                                                            \begin{array}{c}
                                                                                                              0 \\
                                                                                                              1 \\
                                                                                                            \end{array}
                                                                                                          \right)
 $ column at the middle in $\overline{M_0^{2i}}$, then those actions on $\overline{M_0^{2i}}$  can be seen as the actions in the case $g=4k+1$.
 If we add a $\left(
                                                                                                            \begin{array}{c}
                                                                                                              0 \\
                                                                                                              0 \\
                                                                                                            \end{array}
                                                                                                          \right)
 $ column at the middle in $\overline{M_0^{2i+1}}$, then those actions on $\overline{M_0^{2i+1}}$  can be seen as the actions in the case $g=4k+1$.}

 Thus by induction,
 \begin{align*}
&\overline{M_0^{g-1}}\circ \sigma_{g}\sigma_{g+1}\sigma_{g+2}\sigma_{g+1}\sigma_g \\
=&\overline{M_0^{g-1}}\circ \sigma_{g+1}\cdot\sigma_{g+2}\sigma_{g} \\
=&\left(\begin{array}{ccccc}
  A_{k}  & \begin{array}{cc}
                0 & 0 \\
                0 & 0
              \end{array}
 & A_{k}
\end{array}\right) \\
=& \overline{M_0^g},
\end{align*}
then we get $\overline{M_0}\circ\tau=\overline{M_0^g}\circ \sigma_{g+1}=\overline{M_0}$.

Now let $g=4k+2$. It is easy to see that the transpositions $\{\sigma_i | 1 \leq i \leq 2g+1, i \neq 4k+3+2m=g+1+2m\}$ generate a subgroup of the isotropy group at $\overline{M_m}$,
$0 \leq m \leq 2k+1$.

We consider the $\tau$ action on $$\overline{M_0}=\left(\begin{array}{ccc}
A_{k+1}
  & A_{k+1}\end{array}\right).$$

  Since transpositions $$(i,2g+3-i)=\sigma_{i}\sigma_{i+1}\cdots\sigma_{2g+1-i}\sigma_{2g+2-i}\sigma_{2g+1-i}\cdots\sigma_{i+1}\sigma_i,$$ $\tau$ can be represented by the product of transpositions
\begin{align*}
\tau= &\sigma_1\sigma_2\cdots\sigma_{2g+1}\cdots\sigma_2\sigma_1 \\
& \cdot\sigma_2\sigma_3\cdots\sigma_{2g}\cdots\sigma_3\sigma_2 \\
& \cdots \\
& \cdot\sigma_{g}\sigma_{g+1}\sigma_{g+2}\sigma_{g+1}\sigma_{g}\\
& \cdot\sigma_{g+1}.
 \end{align*}

 First we consider the $\sigma_1\sigma_2\cdots\sigma_{2g+1}\cdots\sigma_2\sigma_1$ action on $\overline{M_0}$,
\begin{align*}
&\overline{M_0}\circ \sigma_1\sigma_2\cdots\sigma_{2g+1}\cdots\sigma_2\sigma_1 \\
=&\overline{M_0}\circ \sigma_{g+1}\cdot\sigma_{g+2}\sigma_{g+3}\cdots\sigma_{2g+1}\sigma_{g}\sigma_{g-1}\cdots\sigma_1 \\
=&\left(
     \begin{array}{ccccc}
         \begin{array}{c}
           1 \\
           0 \\
         \end{array}
        & A_{k} &
                      \begin{array}{ccc}
                        1 & 1 \\
                        0 & 0 \\
                      \end{array}
         & A_{k} & \begin{array}{c}
           1 \\
           0 \\
         \end{array}
     \end{array}
   \right) \\
=& \overline{M_0^1}.
\end{align*}

Second we cosider the $\sigma_2\sigma_3\cdots\sigma_{2g}\cdots\sigma_3\sigma_2$ action on $\overline{M_0^1}$,
\begin{align*}
&\overline{M_0^1}\circ \sigma_2\sigma_3\cdots\sigma_{2g}\cdots\sigma_3\sigma_2 \\
=&\overline{M_0^1}\circ \sigma_{g+1}\cdot\sigma_{g+2}\sigma_{g+3}\cdots\sigma_{2g}\sigma_{g}\sigma_{g-1}\cdots\sigma_2 \\
=&\left(
     \begin{array}{ccccc}
         \begin{array}{cc}
          0 & 1 \\
          1 & 0 \\
         \end{array}
        & A_{k} & A_{k} & \begin{array}{cc}
           1 & 0 \\
           0 & 1 \\
         \end{array}
     \end{array}
   \right) \\
=& \overline{M_0^2} \\
& \cdots
\end{align*}

{ In fact, if we add  a  $\left(
                                                                                                            \begin{array}{c}
                                                                                                              0 \\
                                                                                                              0 \\
                                                                                                            \end{array}
                                                                                                          \right)
 $ column at the middle in $\overline{M_0^{2i}}$, then those actions on $\overline{M_0^{2i}}$  can be seen as the actions in the case $g=4(k+1)-1=4k+3$.
 If we add a $\left(
                                                                                                            \begin{array}{c}
                                                                                                              0 \\
                                                                                                              1 \\
                                                                                                            \end{array}
                                                                                                          \right)
 $ column at the middle in $\overline{M_0^{2i+1}}$, then those actions on $\overline{M_0^{2i+1}}$  can be seen as the actions in the case $g=4k+3$.}

 Thus by induction,
 \begin{align*}
&\overline{M_0^{g-1}}\circ \sigma_{g}\sigma_{g+1}\sigma_{g+2}\sigma_{g+1}\sigma_g \\
=&\overline{M_0^{g-1}}\circ \sigma_{g+1}\cdot\sigma_{g+2}\sigma_{g} \\
=&\left(\begin{array}{ccccc}
  A_{k}  & \begin{array}{cccc}
               1 & 0 & 0 & 1 \\
               0 & 1 & 1 & 0
              \end{array}
 & A_{k}
\end{array}\right) \\
=& \overline{M_0^g},
\end{align*}
then we get $\overline{M_0}\circ\tau=\overline{M_0^g}\circ \sigma_{g+1}=\overline{M_0}$.

The theorem follows. \qb

In conclusion, for both of the cases $g=2k-1$ and $g=2k$, $k \geq 2$,  there exists an element $\sigma \in \mathfrak{S}_{2g+2}$ such that $M_m \circ \sigma=\overline{M_m}$, since $M_m$ and $\overline{M_m}$ are in the same orbit $\mathcal{C}_m$.
Thus there is a subgroup of the isotopy group at $M_m$ which can be represented by $\sigma \cdot \mathfrak{G}_m \cdot \sigma^{-1}$.

Now we will finish our proof of Theorem $3.5$ and Theorem $3.6$.

{\bf Proof of Theorem 3.5 and Theorem 3.6:}
 {Denote by $|S|$ the number of elements of a set $S$.  Since any spin structure $c$ uniquely depends on the choice of $f \in H^1(\Sigma_g,\mathbb{Z}/2)$, we have $|\mathcal{C}|=2^{2g}$.}

Let $g=2k-1$. From Theorem $4.2$, we have $$|\mathcal{C}|= 2^{2g}\leq|\mathcal{C}_0|+|\mathcal{C}_1|+\cdots +|\mathcal{C}_k|,$$
and from Theorem $4.4$, we have $$|\mathfrak{G}_m|= (g+1+2m)!\cdot(g+1-2m)!\text{, } |\mathfrak{G}_0|= 2\cdot(g+1)!\cdot(g+1)!,$$
$$|\mathcal{C}_m| \leq |\mathfrak{S}_{2g+2}|/|\mathfrak{G}_m|=\left(
                                                                                                  \begin{array}{c}
                                                                                                    2g+2 \\
                                                                                                    g+1-2m \\
                                                                                                  \end{array}
                                                                                                \right)
$$ for $1 \leq m \leq k$, and $$|\mathcal{C}_0| \leq |\mathfrak{S}_{2g+2}|/|\mathfrak{G}_0|=\frac{1}{2}\left(
                                                                                                  \begin{array}{c}
                                                                                                    2g+2 \\
                                                                                                    g+1 \\
                                                                                                  \end{array}
                                                                                                \right).$$ Here $\left(
                                                                                                  \begin{array}{c}
                                                                                                    n \\
                                                                                                    i \\
                                                                                                  \end{array}
                                                                                                \right)=\frac{(n)!}{(i)!(n-i)!}$ is the binomial coefficient.

Thus we have
\begin{align*}
|\mathcal{C}_0|+|\mathcal{C}_1|+\cdots +|\mathcal{C}_k|\leq & \frac{1}{2}\left(
                                                                                                  \begin{array}{c}
                                                                                                    2g+2 \\
                                                                                                    g+1 \\
                                                                                                  \end{array}
                                                                                                \right) + \left(
                                                                                                  \begin{array}{c}
                                                                                                    2g+2 \\
                                                                                                    g-1 \\
                                                                                                  \end{array}
                                                                                                \right) +\cdots +\left(
                                                                                                  \begin{array}{c}
                                                                                                    2g+2 \\
                                                                                                    0 \\
                                                                                                  \end{array}
                                                                                                \right)\\
=& \left(
                                                                                                  \begin{array}{c}
                                                                                                    2g+1 \\
                                                                                                    g \\
                                                                                                  \end{array}
                                                                                                \right)+(\left(
                                                                                                  \begin{array}{c}
                                                                                                    2g+1 \\
                                                                                                    g-1 \\
                                                                                                  \end{array}
                                                                                                \right)+\left(
                                                                                                  \begin{array}{c}
                                                                                                    2g+1 \\
                                                                                                    g-2 \\
                                                                                                  \end{array}
                                                                                                \right))+\cdots \\
& +(\left(
                                                                                                  \begin{array}{c}
                                                                                                    2g+1 \\
                                                                                                    2 \\
                                                                                                  \end{array}
                                                                                                \right)+\left(
                                                                                                  \begin{array}{c}
                                                                                                    2g+1 \\
                                                                                                    1 \\
                                                                                                  \end{array}
                                                                                                \right))+\left(
                                                                                                  \begin{array}{c}
                                                                                                    2g+2 \\
                                                                                                    0 \\
                                                                                                  \end{array}
                                                                                                \right)\\
=& \frac{1}{2}\sum_{i=0}^{2g+1}\left(
                                 \begin{array}{c}
                                   2g+1 \\
                                   i \\
                                 \end{array}
                               \right)\\
=& 2^{2g}.
\end{align*}
Hence $|\mathcal{C}_0|+|\mathcal{C}_1|+\cdots +|\mathcal{C}_k|=2^{2g}=|\mathcal{C}|$, Theorem $3.5$ follows.

Now let $g=2k$. Similarly, from Theorem $4.2$, we have $$|\mathcal{C}|= 2^{2g}\leq|\mathcal{C}_0|+|\mathcal{C}_1|+\cdots +|\mathcal{C}_k|,$$
and from Theorem $4.5$, we have $$|\mathfrak{G}_m|= (g+1+2m)!\cdot(g+1-2m)!\text{, } |\mathfrak{G}_0|= 2\cdot(g+1)!\cdot(g+1)!,$$
$$|\mathcal{C}_m| \leq |\mathfrak{S}_{2g+2}|/|\mathfrak{G}_m|=\left(
                                                                                                  \begin{array}{c}
                                                                                                    2g+2 \\
                                                                                                    g+1-2m \\
                                                                                                  \end{array}
                                                                                                \right)
$$ for $1 \leq m \leq k$, and $$|\mathcal{C}_0| \leq |\mathfrak{S}_{2g+2}|/|\mathfrak{G}_0|=\frac{1}{2}\left(
                                                                                                  \begin{array}{c}
                                                                                                    2g+2 \\
                                                                                                    g+1 \\
                                                                                                  \end{array}
                                                                                                \right).$$

Thus we have
\begin{align*}
|\mathcal{C}_0|+|\mathcal{C}_1|+\cdots +|\mathcal{C}_k|\leq & \frac{1}{2}\left(
                                                                                                  \begin{array}{c}
                                                                                                    2g+2 \\
                                                                                                    g+1 \\
                                                                                                  \end{array}
                                                                                                \right) + \left(
                                                                                                  \begin{array}{c}
                                                                                                    2g+2 \\
                                                                                                    g-1 \\
                                                                                                  \end{array}
                                                                                                \right) +\cdots +\left(
                                                                                                  \begin{array}{c}
                                                                                                    2g+2 \\
                                                                                                    1 \\
                                                                                                  \end{array}
                                                                                                \right)\\
=& \left(
                                                                                                  \begin{array}{c}
                                                                                                    2g+1 \\
                                                                                                    g \\
                                                                                                  \end{array}
                                                                                                \right)+(\left(
                                                                                                  \begin{array}{c}
                                                                                                    2g+1 \\
                                                                                                    g-1 \\
                                                                                                  \end{array}
                                                                                                \right)+\left(
                                                                                                  \begin{array}{c}
                                                                                                    2g+1 \\
                                                                                                    g-2 \\
                                                                                                  \end{array}
                                                                                                \right))+\cdots \\
& +(                                                                                            \left(
                                                                                                  \begin{array}{c}
                                                                                                    2g+1 \\
                                                                                                    1 \\
                                                                                                  \end{array}
                                                                                                \right)+\left(
                                                                                                  \begin{array}{c}
                                                                                                    2g+1 \\
                                                                                                    0 \\
                                                                                                  \end{array}
                                                                                                \right))\\
=& \frac{1}{2}\sum_{i=0}^{2g+1}\left(
                                 \begin{array}{c}
                                   2g+1 \\
                                   i \\
                                 \end{array}
                               \right)\\
=& 2^{2g}.
\end{align*}
Hence $|\mathcal{C}_0|+|\mathcal{C}_1|+\cdots +|\mathcal{C}_k|=2^{2g}=|\mathcal{C}|$, Theorem $3.6$ follows. \qb

Since the number of $\mathfrak{S}_{2g+2}-$orbits of spin structures of $\Sigma_g$ is exactly $k+1$ for $g=2k-1$ or $g=2k$, by using the above inequations, we can determine the isotropy group of each orbit immediately.

{\bf Corollary 4.6.} {\it The subgroups which we mentioned at Theorem $4.4$ and Theorem $4.5$ are exactly the isotropy groups at those spin structures.}


\begin{thebibliography}{10}
\bibitem{A68}
V.~I.~Arnold, A remark on the ramification of hyperelliptic integrals as functions of parameters, \emph{Funct. Anal. Appl.} {\bf 2}(1968), 187--189.



\bibitem{J80}
D.~Johnson, Spin structures and quadratic forms on surfaces, \emph{J. London Math. Soc.} {\bf 22}(1980), no. 2, 365--373.

\bibitem{M84}
D.~Mumford, \emph{Tata Lectures on Theta}, II, Birkh$\ddot{a}$user, Boston-Basel-Stuttgart, 1984.


\end{thebibliography}
\end{document}